\documentstyle[11pt]{amsart}

\theoremstyle{plain} \newtheorem{theorem}{Theorem}[section] \newtheorem{lemma}[theorem]{Lemma}
\newtheorem{proposition}[theorem]{Proposition} \newtheorem{corollaire}[theorem]{Corollary}
\theoremstyle{definition}

\theoremstyle{remark}

\theoremstyle{definition}

\theoremstyle{remark}
\newtheorem{remark}{Remark}

\newtheorem{definition}{Definition}[section]
\def\cal{\mathcal}
\def\C{{\Bbb C}}
\def\N{{\Bbb N}}
\newcommand{\M}{{\cal M}}
\def\K{{\Bbb K}}

\begin{document}
\title{On the convergence of formal mappings\footnote{To appear in Comm. Anal. Geom.}}  
\author{Nordine Mir }
\address{Universit\'e de  Rouen, Laboratoire de Math\'ematiques Rapha\"el Salem, CNRS, Site Colbert 
76821 Mont Saint Aignan France }
\email{Nordine.Mir@@univ-rouen.fr}
\commby{Rouen}
\date{}


\bibliographystyle{plain}

\begin{abstract} Let $f :  (M,p)\rightarrow (M',p')$ be a formal
(holomorphic) nondegenerate map, i.e.  with
formal holomorphic Jacobian $J_f$ not
identically vanishing, between two germs of real
analytic generic submanifolds in
$\C^n$, $n\geq 2$, $p'=f(p)$.  Assuming the target
manifold to be real algebraic, and the source manifold
to be minimal at $p$ in the sense of Tumanov, we prove
the convergence of the so-called reflection mapping
associated to $f$.  From this, we deduce the
convergence of such mappings from minimal  real
analytic generic  submanifolds into real algebraic
holomorphically nondegenerate ones, as well as related
results on partial convergence of such maps.  For the
proofs, we establish a principle of analyticity for
formal CR power series. This principle can be used to
reobtain the convergence of formal mappings of real
analytic CR manifolds under a standard nondegeneracy
condition. 
\end{abstract}

\maketitle

\section{Introduction}\label{sec0}
\setcounter{equation}{0}
A formal (holomorphic) mapping $f : (\C^n,p)\rightarrow (\C^{n'},p')$, $(p,p')\in \C^n \times
\C^{n'}$, $n,n'\geq 1$, is a vector $(f_1,\ldots,f_{n'})$ where each $f_j \in \C[[z-p]]$ is a formal
holomorphic power series in $z-p$, and $f(p)=p'$. In the case $n=n'$, a formal mapping $f$ is called {\it nondegenerate}
if its formal holomorphic Jacobian $J_f$
does not vanish identically as a formal power
series in
$z-p$. An important class of nondegenerate formal maps
$f$ consists of those which are invertible, namely
those for which $J_f(p)\not =0$. We call such maps {\it
formal equivalences} or {\it formal invertible maps}.
If $M,M'$ are two smooth real real analytic generic
submanifolds in $\C^n$ and
$\C^{n'}$ respectively (through
$p,p'$ respectively) and of real codimension $c$ and
$c'$ respectively,  we say that a formal mapping $f :
(\C^n,p)\rightarrow (\C^{n'},p')$ sends $M$ into
$M'$ if 
$$\rho'(f(z),\overline{f(z)})=a(z,\bar{z})\cdot \rho (z,\bar{z}),$$
where $\rho =(\rho_1,\ldots,\rho_c)$ and
$\rho'=(\rho'_1,\ldots,\rho'_{c'})$ are local real
analytic defining functions for $M,M'$ respectively and
$a(z,\bar z)$ is a $c'\times c$ matrix with
entries in
$\C[[z-p,\bar z -\bar p]]$.  It is easy to see
that such a definition is independent of the
choice of defining functions for $M$ and $M'$. 
If $f$ is formal mapping as above sending $M$
into $M'$, we may also say that $f$ is a {\it
formal CR mapping}  from $M$ into $M'$. This is
motivated by the fact that, if, in addition, $f$
is convergent near $p$, then $f$ is a real
analytic CR  mapping from $M$ into $M'$.

A natural question which arises is to give necessary and sufficient conditions
so that any formal equivalence between real analytic generic submanifolds
must be
convergent.  Chern and Moser \cite{CM} gave the first results in this direction by
proving the convergence of formal
equivalences between Levi nondegenerate real analytic hypersurfaces.  Later,
Moser and Webster \cite{MW} showed the analyticity of formal invertible mappings
between certain real analytic surfaces of dimension two in
$\C^2$, but which are not CR.  Other related work was
done by Webster \cite{W2} and Gong \cite{Go}.  (We
also refer the reader to the bibliography given in
\cite{BER4} for further information.)  More recently, Baouendi, Ebenfelt and
Rothschild proved the convergence of formal equivalences between minimal {\it
finitely nondegenerate} real analytic generic submanifolds \cite{BER2, BER3}, as well as
between minimal {\it essentially finite} ones\footnote{See \S \ref{sec2}
for precise definitions.}  \cite{BER4} (other situations are
also treated in \cite{BER3, BER4}).  The conditions of finite nondegeneracy and essential
finiteness are closely related to the notion of {\it holomorphic nondegeneracy}
introduced by Stanton \cite{St}.  Let us recall that a connected real analytic
generic submanifold is {\it holomorphically
nondegenerate} if, near any point $p\in M$, there is
no non-trivial holomorphic vector field, with
holomorphic coefficients, tangent to $M$ near $p$.  
Such submanifolds have the property to be generically
essentially-finite in the sense that, for any such
manifold $M$, there always exist a proper real
analytic subvariety $S\subset M$ (which may be empty)
such that any point $p\in M\setminus S$ is essentially
finite.  Moreover, it was observed in \cite{BER2} that
the condition of holomorphic nondegeneracy is
necessary for the convergence of formal equivalences
between real analytic generic submanifolds.  Thus, to
complete the previous results, one has to treat the
case of the non-essentially finite points of such
holomorphically nondegenerate submanifolds. 

In the one-codimensional case, these non-essentially finite points were
recently treated in \cite{M3} where, in particular, it was shown that
any formal CR equivalence between minimal holomorphically nondegenerate real analytic hypersurfaces 
must be convergent. The goal of this paper is to study the higher-codimensional case. Assuming the target
manifold to be real {\it algebraic} i.e. contained in
a real algebraic subvariety of the same dimension, we
establish a result which gives a description of the
analyticity properties of formal CR nondegenerate maps
from minimal real analytic generic submanifolds of
$\C^n$ into real algebraic ones (Theorem \ref{th1.1}
below).
  As in \cite{M2, M3}, we prove that, given a formal
map
$f : (M,p)\rightarrow (M',p')$ with $J_f\not \equiv
0$, if $M$ is minimal at $p$, then the so-called
associated {\it reflection mapping} ({\it cf.}
\cite{H}) must be convergent. As we shall see ({\it
cf.} \S \ref{secpc}), such a result can be seen as a
result of partial convergence for formal CR
nondegenerate maps. This allows one also to deduce the
convergence of such maps from real analytic minimal
generic submanifolds onto real algebraic
holomorphically nondegenerate ones (Theorem
\ref{cor1.2} below).  We should point out that our
arguments give also a quite simple proof of such a
fact (see Proposition \ref{prop6.2}).  In fact, the
algebraicity of the target manifold allows us to use
certain tools from basic field theory that we
introduced in our previous works
\cite{M1, M2}.
\section{Statement of main results}\label{sec1}
\setcounter{equation}{0}
Let $(M',p')\subset \C^n$ be a germ at $p'$ of a real algebraic generic submanifold of CR dimension
$N$ and of real codimension $c$. This means that there
exists $\rho'
(\zeta,\bar{\zeta})=(\rho'_1(\zeta,\bar{\zeta}),
\ldots,\rho'_c(\zeta,\bar{\zeta}))$ $c$ real
polynomials such that near $p'$
$$M'=\{\zeta \in (\C^n,p'):\rho'(\zeta,\bar{\zeta})=0\},$$
with $\bar{\partial}\rho'_1\wedge \ldots \wedge
\bar{\partial}\rho'_c \not =0$, on $M$. We shall
assume, without loss of generality, that $p'$ is the
origin. Then, for any point $\omega$ close to
$0$, one defines its associated {\it Segre variety} to be the 
$n-c$ dimensional complex submanifold 
\begin{equation}\label{sv}
Q'_{\omega}=\{\zeta \in
(\C^n,0):
\rho'(\zeta,\bar{\omega})=0\}.
\end{equation} Moreover, since $M'$
is generic, renumbering the coordinates if necessary,
and after applying the implicit function theorem, one
can assume that any Segre variety can be described as
a graph of the form
$$Q'_{\omega}=\{\zeta \in (\C^n,0): \bar{\zeta}^*=\bar{\Phi}'(\omega,\bar{\zeta}')\},
\quad \zeta =(\zeta',\zeta^*)\in \C^N\times
\C^c,$$
$\bar{\Phi}'=(\bar{\Phi}'_1,\ldots,\bar{\Phi}'_c)$
denoting a convergent power series mapping near $0\in
\C^{n+N}$, with $\bar{\Phi}'(0)=0$. Our main result is
the following.

\begin{theorem}\label{th1.1} Let $f:(M,0)\rightarrow (M',0)$ be a formal nondegenerate CR map between two 
germs at 0 of real
analytic generic submanifolds in $\C^n$ of the same
CR dimension.  Assume that
$M$ is minimal at $0\in M$ and that $M'$ is real
algebraic. Then, the formal holomorphic map
$$\C^n\times \C^N \ni (z,\theta)\mapsto
\bar{\Phi}'(f(z),\theta) \in \C^c$$ is
convergent. \end{theorem} Such a result was
established in \cite{M3}  in the
one-codimensional case (for unbranched
mappings) without assuming that the target manifold 
$M'$ is real algebraic.
As in \cite{M3}, Theorem \ref{th1.1} allows us
to derive the following convergence result.

\begin{theorem}\label{cor1.2} 
Let $f:(M,0)\rightarrow (M',0)$ be a formal nondegenerate CR map between two 
germs at 0 of real
analytic generic submanifolds in $\C^n$ of the same
CR dimension.  Assume that
$M$ is minimal at $0\in M$ and that $M'$ is real
algebraic and holomorphically nondegenerate. Then,
$f$ is convergent. 
\end{theorem} As mentioned in the introduction,
Theorem \ref{cor1.2}, for unbranched maps, follows
from \cite{M3} in the hypersurface case, but in the
higher codimensional case the result is new and was
not previously known even in the case
 where $M$ and $M'$ are both algebraic. Another application of Theorem \ref{th1.1} 
 is given in \S \ref{secpc} and deals with partial convergence of formal CR nondegenerate maps. For this, we
  refer the reader to Theorem \ref{thpc0} and Corollary
\ref{thpc}.

Our approach for proving Theorem \ref{th1.1} is
essentially based  on two steps. The first step
is a formulation of the reflection principle via the jet method and follows \cite{M2}.  The general idea is to
show that, under the assumptions of Theorem
\ref{th1.1}, the composition of any component of the
Segre variety map of $M'$ (as defined in \S
\ref{sec4}) with the map $f$ satisfy certain
polynomial equations restricted on $M$, and more
precisely, is algebraic over a certain field of formal
power series.  The second step consists in proving
that a formal CR power series (i.e. a formal
holomorphic power series) which satisfies such a
polynomial identity, is necessarily convergent (Theorem
\ref{th5.1}).  This is based on the theory of Segre
sets by Baouendi, Ebenfelt and Rothschild \cite{BER1},
and on some of their techniques of propagation.  One
should, however, point out several differences with
the methods of
\cite{BER1, BER4} (see especially Proposition
\ref{prop5.3}).  \par \smallskip The paper is
organized as follows.  \S \ref{sec2} contain some
background material for the reader's convenience.  In
\S
\ref{sec4}, we use some ideas from \cite{M2} to prove
our reflection identities.  \S \ref{sec5} is devoted to
the proof of a principle of analyticity  for formal CR
functions. Such a result (Theorem \ref{th5.1}) seems
to us interesting in itself. In
\S
\ref{sec6}, we prove the main results of the paper. In
\S
\ref{sec+}, we formulate some remarks concerning
Theorem
\ref{cor1.2} which show that, under the assumptions of
that Theorem, the convergence of formal nondegenerate
maps can be derived in a quite simple manner.  In \S
\ref{secpc}, we apply Theorem
\ref{th1.1} to the study of partial convergence for
formal CR maps. Finally, in
\S \ref{sec7}, we apply the principle proved in \S
\ref{sec5} to establish the convergence of formal
mappings between real analytic CR manifolds under a
standard nondegeneracy condition.

{\sc Acknowledgements.}  This work has been
completed while I was invited by the department of
Mathematics of the university of Wuppertal, Germany,
during the period May-July 1999.  I would like to
thank Prof.  K.  Diederich, C.  Eppel and Prof.  G. 
Herbort for arranging my visit.  I would like also
to thank Prof.  K.  Diederich for interesting
conversations.  I am indebted to Prof. M.  Derridj
for his precious help, all his encouragements and
for having spent many of his time thinking of my
numerous questions.  Finally, I wish also to thank the
referee for many valuable comments and helpful
suggestions.

\section{Preliminaries, notations and definitions.}\label{sec2} \setcounter{equation}{0} \subsection{Finite
nondegeneracy, essential finiteness and holomorphic nondegeneracy of real
analytic generic submanifolds.}\label{ssec2.1} Let $M$ be a real analytic generic submanifold through $p\in \C^n$,
of CR dimension $N$ and of real codimension $c$.  We
shall always assume that $c,N\geq 1$, and, for
convenience, that the reference point $p$ is the
origin.  Let $\rho =(\rho_1,\ldots,\rho_c)$ be a set
of real analytic defining functions for $M$ near 0,
i.e.  \begin{equation}\label{eq2.1} M =\{z\in
(\C^n,0):\rho (z,\bar{z})=0\}, \end{equation} with
$\bar{\partial} \rho_1 \wedge \ldots \wedge
\bar{\partial} \rho_c \not =0,\ {\rm on}\ M$. 
The complexification ${\cal M}$ of $M$ is the
$2n-c$-dimensional complex submanifold of
$\C^{2n}$ given by 
\begin{equation}\label{eqcomp}
\M =\{(z,w)\in (\C^{2n},0):\rho (z,w)=0\}.
\end{equation}
We shall assume, without loss of generality,
that the matrix $\partial \rho/\partial z^*$ is not singular at the origin,
$z=(z',z^*)\in \C^{N}\times \C^c$.  In this case, we
define the following vector fields tangent to $\M$, 
\begin{equation}\label{vf2}
{\cal
L}_j=\frac{\partial}{\partial w_j}- \rho_{w_j}(z,w)\left[\frac{\partial \rho}{\partial w^*}(z,w)\right]^{-1}
\frac{\partial }{\partial w^*},\ j=1,\ldots,N,
\end{equation}
which are the complexifications of the $(0,1)$ vector fields tangent to $M$. Let us also recall that the invariant Segre varieties
attached to $M$ are defined by $$Q_w=\{z\in (\C^n,0):\rho (z,\bar{w})=0\},$$ for $w$ close to 0.  A
fundamental map which arises in the mapping problems is the so-called {\it variety Segre map}
$\lambda :  w\mapsto Q_w$ ({\it cf.}  \cite{DW, DF,
DP}). A real analytic generic submanifold $M$ is
called {\it finitely nondegenerate} at $p=0\in M$ if
$${\rm Span}_{\C}\{{\cal L}^{\alpha}\rho_{j,z}(p,\bar{p}):\alpha \in \N^N,\ 1\leq j\leq d\}=\C^n.$$
Here, for $1\leq j\leq d$, $\rho_{j,z}$ denotes the complex gradient of $\rho_j$ with respect to $z$. In this
case, one can show that the Segre variety map
$\lambda$ is actually one-to-one near $p=0$. More
generally, $M$ is called {\it essentially finite} at
$0\in M$ if the Segre variety map $\lambda$ is
finite-to-one near
$0$ \cite{DW, BJT}.  The interest of such
conditions lies in the fact that, given a
holomorphically nondegenerate generic real
analytic submanifold
$M$ (as defined in the introduction), the set of
finitely nondegenerate or essentially finite points is
always, at least, dense in $M$ (see \cite{BER1}). 
Furthermore, the set of points $w\in M$ such that
$\lambda^{-1}(Q_{w})$ is positive-dimensional forms a
proper (possibly empty) real analytic  subvariety
$S\subset M$, provided that the submanifold $M$ is
holomorphically nondegenerate. This set of points is
precisely the set of {\it non-essentially finite}
points of
$M$.  

\subsection{Minimality condition in terms of Segre sets}\label{ssec2.2}
Another nondegeneracy condition which will be
used in this paper is the minimality condition
introduced by Tumanov
\cite{Tu}.  Let us recall that a real analytic generic submanifold $M$ is said to be {\it minimal} at $p\in M$ (or of {\it finite
type} in the sense of Kohn and Bloom-Graham) if there is no proper CR
submanifold contained in $M$ through $p$, and with the same Cauchy-Riemann dimension.  In order to give a
characterization of minimality for real analytic CR manifolds, Baouendi, Ebenfelt and Rothschild introduced
the so-called {\it Segre sets} in \cite{BER1}.  These
sets will play an important role in our proofs.  They
are defined as follows.  Define the first Segre set
$N_1(p)$ attached to $M$ at $p\in M$ to be the
classical Segre variety $Q_p$.  Inductively, for $k\in
\N$, define $$N_{k+1}(p)=\bigcup_{q\in N_k(p)}Q_q.$$
Recall that the sets $N_j(p)$ are in general not
analytic for $j>1$.  If $M$ is given by (\ref{eq2.1})
near $p=0$, as in \S \ref{ssec2.1}, by the implicit
function theorem, one can choose coordinates
$z=(z',z^*)\in \C^N\times \C^{c}$ so that any
Segre variety can be described as a graph as
follows
$$Q_w=\{z\in (\C^n,0):z^*=\Phi (\bar{w},z')\},$$
$\Phi =(\Phi_1,\ldots,\Phi_c)$ being a
$\C^c$-valued holomorphic map near $0\in
\C^{N+n}$,
$\Phi (0)=0$. 
The reality of $M$ also implies that
\begin{equation}\label{eq2000} \Phi (w',
\bar{\Phi}(z', z^*,w'),z')\equiv z^*,\ (z,w')\in
\C^n\times \C^N.  \end{equation} (Here and in
what follows, for any formal power series
$g(x)\in \C[[x]]$,
$x=(x_1,\ldots,x_k)$, $\bar{g}(x)$ is the formal power series obtained by taking the complex
conjugates of the coefficients of $g$.)  The
coordinates are said to be {\it normal} for $M$ if,
moreover, the condition $\Phi (z,0)\equiv z^*$
holds.  It is well known ({\it cf.}  \cite{CM,
BJT}) that given a real analytic generic
submanifold $M$, one can always find such
coordinates.  With these notations (and without
assuming that the $z$-coordinates are normal for
$M$), the Segre sets can be parametrized by the
following mappings $(v_k)_{k\in \N}$, called the
{\it Segre sets mappings}.  First set $v_0:=0\in \C^n$. Inductively,
$N_{2k+1}$,
$k\geq 0$, can be parametrized by the map
\begin{multline}\label{eqss1}
(\C^{(2k+1)N},0)\ni (t_1,t_2,\ldots,t_{2k+1})
\mapsto \\v_{2k+1}(t_1,\ldots,t_{2k+1}):=\left( t_1,
\Phi (
\bar{v}_{2k}(t_2,\ldots,t_{2k+1}),t_1)\right)
\end{multline} and $N_{2k}$ by
\begin{multline}\label{eqss2} (\C^{2kN},0)\ni
(t_1,t_2,\ldots,t_{2k}) \mapsto
\\v_{2k}(t_1,\ldots,t_{2k}):= (t_1,\Phi
(\bar{v}_{2k-1}(t_2,\ldots,t_{2k}),t_1)).
\end{multline} Notice that, for any nonnegative integer
$b$,
$(v_{b+1}(t_1,\ldots,t_{b+1}),\bar{v}_{b}(t_2,\ldots,t_{b+1}))\in {\cal M}$.  We can now state a useful
characterization of minimality which is contained in \cite{BERb}.

\begin{theorem}\label{th2.1} {\rm \cite{BERb}} If $M$ is minimal at 0, there exists $d_0\in \N$ large
enough such that, in any neighborhood ${\cal O}$ of $0\in \C^{d_0N}$, there exists
$(t_1^0,\ldots,t_{d_0}^0)\in {\cal O}$, such that $v_{d_0}(t_1^{0},\ldots,t_{d_0}^0)=0$ and such that
$v_{d_0}$ is submersive at $(t_1^0,\ldots,t_{d_0}^0)$.  \end{theorem}

\section{Formal nondegenerate CR maps with
values in real algebraic CR manifolds}\label{sec4}
\setcounter{equation}{0} 
\subsection{Real algebraic CR manifolds and field
extensions}\label{ssec4.1} In this section, we collect
and recall  some
facts from \cite{MT, M1, M2} which will be used in
the proof of Theorem \ref{th1.1}.

As in \S \ref{sec1}, let $(M',p')$ be a germ through
$p'=0\in \C^n$ of a real algebraic generic submanifold
of
CR dimension $N$ and of real codimension $c$.
Following the notations of that section, let
$\rho'=(\rho'_1,\ldots,\rho'_c)$ be a set of  
defining real polynomials for $M'$ near 0.  Thus, 
\begin{equation}\label{defM'}
M' =\{\zeta\in
(\C^n,0):\rho' (\zeta,\bar{\zeta})=0\},
\end{equation} with
$\bar{\partial}
\rho'_1 \wedge \ldots \wedge \bar{\partial} \rho'_c
\not =0,\ {\rm on}\ M'$.  We can assume that the
coordinates $\zeta$ at the target space are chosen so
that if $\zeta=(\zeta',\zeta^*)\in \C^N\times
\C^{c}$, the matrix $\partial \rho'/\partial \zeta^*$
is not singular at the origin.  This allows one to
represent
$M'$ as follows  
\begin{equation}\label{eqM"}
M'=\{\zeta\in
(\C^n,0):\bar{\zeta}^*=\bar{\Phi}'(\zeta,\bar{\zeta}')\},
\end{equation} 
$\bar{\Phi}'
=(\bar{\Phi}'_1,\ldots,\bar{\Phi}'_c)$ being a
$\C^c$-valued holomorphic {\it algebraic} map
near $0\in \C^{N+n}$ with
$\bar{\Phi}'(0)=0$. (We recall here that a holomorphic
function in $k$ variables near 0 is called algebraic
if it is algebraic over the quotient field of the
polynomials in $k$ indeterminates.) Write, for
$\nu=1,\ldots,c$, the expansion 
\begin{equation}\label{taylor}
\bar{\Phi}'_{\nu}
(\omega,\theta)=\sum_{\beta
\in \N^{N}}q_{\beta,\nu}(\omega)
\theta^{\beta}.
\end{equation}
Here, $\omega \in \C^n$ and $\theta \in \C^N$. We also
write
\begin{equation}\label{note}
\bar{\Phi}'_{\theta^{\alpha}}(\omega,\theta)=
(\bar{\Phi}'_{\theta^{\alpha},1}(\omega,\theta),\ldots,\bar{\Phi}'_{\theta^{\alpha},c}(\omega,\theta))=(\partial_{\theta}^{\alpha}\bar
\Phi'_1(\omega,\theta),\ldots,\partial_{\theta}^{\alpha}\bar\Phi'_c(\omega,\theta)).
\end{equation} With these notations, the  Segre variety
map
$\lambda':(\C^n,0) \ni \omega \mapsto
Q'_{\omega}$ associated to $M'$ can be identified 
with the holomorphic map
\begin{equation}\label{eqSVM}
(\C^n,0) \ni \omega \mapsto
\left(q_{\beta,\nu}(\omega)\right)_{\beta \in \N^N
\atop 1\leq \nu \leq c}.
\end{equation} Here, the Segre variety
$Q'_{\omega}$, for $\omega$ close to 0, is defined
by (\ref{sv}). The family of holomorphic
algebraic functions defined by (\ref{eqSVM}) will be
denoted
${\cal C}$.
For $k\in \N^*$, let ${\cal F}_k$ be the quotient field of the germs at $0\in \C^k$ of
algebraic functions in $\C^k$.  For any positive integer $l\in \N$, we define ${\cal P}_l$ to be the smallest
field contained in ${\cal F}_{N+n}$ and containing $\C$ and the family $(\theta,
\bar{\Phi}'_{\theta^{\beta},j}(\omega,\theta))_{j=1,
\ldots,c,
|\beta|\leq l}$.  We then define ${\cal P}\subset
{\cal F}_{N+n}$ to be the set 
\begin{equation}\label{eqP}
{\cal P}=\cup_{l\in
\N}{\cal P}_l.
\end{equation} One can easily check that
${\cal P}$ is also a subfield of $ {\cal F}_{N+n}$, since, for any $l$, ${\cal P}_l\subset {\cal P}_{l+1}$.
By definition, an element $b=b(\omega,\theta)\in
{\cal F}_{N+n}$ belongs to ${\cal P}$ if there
exists a positive integer $l$ and two holomorphic
polynomials $Q_1$ and $Q_2$ such that
$Q_2\left((\bar{\Phi}'_{\theta^{\beta},j}(\omega,\theta))_
{j\leq
c, |\beta|\leq l},\theta \right)\not \equiv 0$ in
${\cal F}_{N+n}$ and such that
$$b=b(\omega,\theta)=\displaystyle
\frac{Q_1\left((\bar{\Phi}'_{\theta^{\beta},j}
(\omega,\theta))_{j\leq
c, |\beta|\leq
l},\theta
\right)}{Q_2\left((\bar{\Phi}'_{\theta^{\beta},j}(\omega,\theta))_{j\leq
c, |\beta|\leq l},\theta \right)}.$$ We need to state
the following proposition, established in \cite{M2}
(Proposition 1) in the hypersurface case, but which
follows with the same proof in the higher
codimensional case.

\begin{proposition}\label{prop4.1} Let $M'$ be a
real algebraic generic submanifold of CR
dimension $N$ through the origin in $\C^n$. Assume
that
$M'$ is given near 0 by
$(\ref{eqM"})$. Let
${\cal C}$ be the family of algebraic holomorphic
functions (in
$n$ variables) defined by $(\ref{eqSVM})$ and
${\cal P}$ be the field of algebraic holomorphic
functions (in $N+n$ variables) defined by
$(\ref{eqP})$. Then, the following holds. The family
${\cal C}$ is contained in the algebraic closure
of ${\cal P}$, and hence, the algebraic closure of
${\cal C}$ is contained in the algebraic closure
of ${\cal P}$.
\end{proposition} \begin{remark}\label{rk1} {\rm An
inspection of the proof of Proposition 1 from
\cite{M2} shows that there exists $l_0$, which depends
only on $M'$ such that
${\cal C}$ is contained in the algebraic closure of
${\cal P}_{l_0}$.  Moreover, if $M'$ is
holomorphically nondegenerate, $l_0$ is nothing else
than the so-called {\it Levi-type} of $M'$
(see \cite{BER1}).  Indeed, we define $l_0$ as follows.  Consider, for any positive integer $l$,
the map
$\psi_l: (\C^{N+n},0)\ni (\omega,\theta)\mapsto
(\theta,(\bar{\Phi}'_{\theta^{\beta},j}(\omega,\theta))_{j\leq
c,|\beta|\leq l})$ and denote by $r_l$ the  generic
rank of such a map.  Finally, put $r(M')={\rm
max}_{l\in \N}\ r_l$.  Then, $l_0={\rm inf}\{l\in
\N:r_l=r(M')\}$.  When $M'$ is holomorphically
nondegenerate, then it is well-known ({\it cf.} 
\cite{BER1}) that, in that case, the integer
$r(M')$ equals $N+n$ and, by definition, $l_0$ is
the Levi-type of $M'$.} \end{remark} 
We recall also the following criterion of holomorphic
nondegeneracy from \cite{MT,M1}.

\begin{theorem}\label{thmir}
Let $M'$ be a real algebraic generic submanifold
through the origin in $\C^n$. Assume also that $M'$
is given near 0 by $(\ref{eqM"})$. Let ${\cal C}$
 be the family of algebraic holomorphic
functions (in
$n$ variables) defined by $(\ref{eqSVM})$ and
${\cal P}$ be the field of algebraic holomorphic
functions (in $N+n$ variables) defined by
$(\ref{eqP})$. Then, the following conditions are
equivalent:
\begin{enumerate}
\item[(i)] $M'$ is holomorphically nondegenerate (at 0)
\item[(ii)] The algebraic closure of the field
${\cal P}$ is ${\cal F}_{N+n}$
\item[(iii)] The algebraic closure of the field
generated by ${\cal C}$ is ${\cal F}_n$
\end{enumerate}
\end{theorem} 
\subsection{Jets and the reflection
principle}\label{ssec4.2} In this section, we
assume that we are in the following setting. Let
$f:  (M,p)\rightarrow (M',p')$ be a formal CR map
between two real analytic generic submanifolds in
$\C^n$, with the same CR dimension $N$ and same real
codimension
$c$.  We assume that
$f$ is a nondegenerate map, i.e. that its formal holomorphic Jacobian 
$J_f$ is not identically vanishing. We 
also assume that $M'$ is a real algebraic generic
submanifold and, without loss of generality, that $p$
and $p'$ are the origin. We use the notations
introduced in \S \ref{sec2} for $M$, and those
introduced in \S \ref{ssec4.1} for $M'$. The goal of
this section is to prove the following proposition.

\begin{proposition}\label{propr3} Let $M\subset \C^n$ be a real
analytic generic submanifold through the origin and
$M'\subset \C^n$ be a real algebraic generic
submanifold through the origin with the same CR
dimension. Assume that $M'$ is given near 0 by
$(\ref{eqM"})$. Let
${\cal C}$ be the family of algebraic functions (in
$n$ variables) associated to $M'$ defined by
$(\ref{eqSVM})$. Let
$\chi
\in {\cal C}$ and
$f:M\rightarrow M'$ be a formal CR map between $M$ and
$M'$ with $J_f\not \equiv 0$.  Then, there exists
$l_0\in \N^*$ (depending only on $M'$), a positive
integer $k_0$ (depending only on $M'$ and $\chi$) and
a family of convergent power series
$\delta_i=\delta_i
\left((\Lambda_{\gamma})_{|\gamma|\leq
l_0},z,w\right)\in
\C\{(\Lambda_{\gamma}-\partial^{\gamma}\bar{f}(0))_{|\gamma|\leq
l_0},z,w\}$, $i=0,\ldots,k_0$, such that the formal
power series identity
$$\sum_{i=0}^{k_0}\delta_i\left((\partial^{\gamma}\bar{f}(w))_{|\gamma|\leq
l_0},z,w\right)\left(\chi\circ f(z)\right)^i =0,$$
holds for
$(z,w)\in {\cal M}$ such that
$\delta_{k_0}\left((\partial^{\gamma}\bar{f}(w))_{|\gamma|\leq
l_0},z,w\right)\not \equiv 0$ in ${\cal M}$.  Here,
${\cal M}$ is the complexification of $M$ as defined
by {\rm (\ref{eqcomp})}.
\end{proposition}
To prove Proposition \ref{propr3}, we will use an approach which is contained in \cite{M2}. We shall first state
several preliminary results needed for its proof. 

 Recall first that the coordinates at the target
space are denoted by $\zeta$. We write
\begin{equation}\label{split}
f=(f',f^*)=(f'_1,\ldots,f'_N,f^*)
\end{equation}
in the
$\zeta=(\zeta',\zeta^*)\in \C^N\times \C^c$
coordinates. Since
$f$ maps formally $M$ into $M'$, there exists
$a(z,\bar{z})$ a $c\times c$ matrix with coefficients
in
$\C[[z,\bar{z}]]$ such that the following formal vectorial identity \begin{equation}\label{eq4.0}
\overline{f^*(z)}-\bar{\Phi}'(f(z),\overline{f'(z)})
=a(z,\bar{z})\cdot \rho (z,\bar{z}), \end{equation}
holds.  Equivalently this gives
\begin{equation}\label{eq4.50} \bar{f}^*(w)-\bar{\Phi}'(f(z),\bar{f}'(w)) =a(z,w)\cdot \rho (z,w),\ {\rm in}\
\C[[z,w]]. \end{equation} 
Define 
\begin{equation}\label{det}
D(z,w)={\rm det}\left({\cal
L}_j\bar{f}'_i(w)\right)_{1\leq i,j\leq N}\in
\C[[z,w]].
\end{equation} Here, ${\cal L}_j$, for $j=1,\ldots,N$,
is the vector field defined by (\ref{vf2}). By applying the vector fields
${\cal L}_j$, $j=1,\ldots,N$, to (\ref{eq4.50}) and Cramer's
rule, one obtains the following known lemma
({\it cf.}
\cite{BER3, DF}).

\begin{lemma}\label{lem4.1} Let $f:(M,0)\rightarrow
(M',0)$ be a formal CR mapping as in Proposition
$\ref{propr3}$. With the notations introduced in
$(\ref{det})$ and $(\ref{note})$, the following holds.
For any multi-index
$\alpha \in \N^{N}$, one has the following
$c$-dimensional formal identity $$
D^{2|\alpha|-1}(z,w)\
\bar{\Phi}'_{\theta^{\alpha}}(f(z),\bar{f}'(w))=V_{\alpha}
\left((\partial^{\beta}\bar{f}(w))_{|\beta|\leq
|\alpha|},z,w\right),$$ 
for $(z,w)\in {\cal M}$. Here,
$V_{\alpha}=(V_{\alpha}^1,\ldots,V_{\alpha}^c)\in
(\C\{(\Lambda_{\beta}-\partial^{\beta}\bar{f}(0))_{|\beta|\leq
|\alpha|},z,w\})^{c}$.  \end{lemma} The following 
lemma contains two known and easy facts. 
 \begin{lemma}\label{lem555} Let $f:(M,0)\rightarrow
(M',0)$ be a formal CR mapping as in Proposition
$\ref{propr3}$. Let $D$ be as in $(\ref{det})$. Then,
the following holds.
\begin{enumerate}
\item[(i)] There exists a convergent power series
$U=U(z,w,(\Lambda_{\beta})_{|\beta|=1})\in
\C\{z,w,(\Lambda_{\beta}-\partial^{\beta}\bar{f}(0))
_{|\beta|= 1})\}$ such that
$D(z,w)=U\left(z,w,(\partial^{\beta}\bar{f}(w))_{|\beta|=
1}\right)$.
\item[(ii)]  $D(z,w)\not =0$ for $(z,w)\in
{\cal M}$. 
\end{enumerate}
\end{lemma} 
{\it Proof of Proposition $\ref{propr3}$.} 
Let ${\cal P}$ be the subfield of
${\cal F}_{N+n}$ defined by (\ref{eqP}). We also
recall that for the positive integer $l_0$
mentioned in Remark \ref{rk1}, ${\cal P}_{l_0}$ is the
smallest field contained in ${\cal F}_{N+n}$ and
containing $\C$ and the family
$(\theta,\bar{\Phi}'_{\theta^{\beta},j}(\omega,\theta))_{1\leq 
j\leq c, |\beta|\leq l_0})$. Let $\chi \in {\cal C}$. Since, by Proposition
\ref{prop4.1},
$\chi$ is algebraic over ${\cal P}$ and, according to 
Remark \ref{rk1}, also over ${\cal P}_{l_0}$, we
obtain the existence of a positive integer $k_0$ and a
family
$(b_j(\omega,\theta))_{0\leq j\leq k_0-1}\subset
{\cal P}_{l_0}$ such that the following identity
\begin{equation}\label{eq2010}
\left(\chi(\omega)\right)^{k_0}+\sum_{j=0}^{k_0-1}
b_j(\omega,\theta)\left(\chi
(\zeta)\right)^j\equiv 0 \end{equation} holds in the
field ${\cal F}_{N+n}$.  By definition, for
$j=0,\ldots,k_0-1$, there exist holomorphic
polynomials
$Q_{1,j}$, $Q_{2,j}$ such that
\begin{equation}\label{eq2011}
Q_{2,j}\left((\bar{\Phi}'_{\theta^{\beta},j}(\omega,\theta))_{j\leq
c, |\beta|\leq l_0},\theta \right)\not \equiv 0,
\end{equation} and such that \begin{equation}\label{eq2012}
b_{j}(\omega,
\theta)=\frac{Q_{1,j}\left((\bar{\Phi}'_{\theta^{\beta},j}(\omega,\theta))_{j\leq
c, |\beta|\leq
l_0},\theta
\right)}{Q_{2,j}\left((\bar{\Phi}'_{\theta^{\beta},j}(\omega,\theta))_{j\leq
c, |\beta|\leq l_0},\theta \right)}.
\end{equation} Now, one sees that (\ref{eq2010}), 
(\ref{eq2011}) and (\ref{eq2012}) imply that there
exist holomorphic polynomials $s_j$,
$j=0,\ldots,k_0$, such that, in some neighborhood of
$0\in \C^{2n-c}$, the following identity
\begin{equation}\label{eqr1}
\sum_{i=0}^{k_0}s_i\left(
(\bar{\Phi}'_{\theta^{\alpha},\mu}(\omega,\theta))_{\mu\leq
c,|\alpha|\leq l_0},\theta \right)
\left(\chi(\omega)\right)^i\equiv 0
\end{equation} holds, with the additional
non-degeneracy condition
\begin{equation}\label{eq4.2} s_{k_0}\left(
(\bar{\Phi}'_{\theta^{\alpha},\mu}(\omega,\theta))_{\mu
\leq c,|\alpha|\leq l_0},\theta \right)\not \equiv
0.  \end{equation} Note that $k_0$ and the family
$(s_i)_{i\leq k_0}$ depend only on $\chi$ and
$M'$. Putting, for $(z,w)\in {\cal M}$,
$\omega=f(z)$ and $\theta=\bar{f}'(w)$ in
(\ref{eqr1}), one obtains the following formal
identity ({\it cf.} 
\cite{M2})
\begin{equation}\label{eq1973} \sum_{i=0}^{k_0}s_i\left(
(\bar{\Phi}'_{\theta^{\alpha},\mu}(f(z),
\bar{f}'(w)))_ {\mu \leq c, |\alpha|\leq
l_0},\bar{f}'(w) \right)
\left((\chi \circ f)(z)\right)^i \equiv 0. 
\end{equation} From Lemma \ref{lem4.1} and Lemma
\ref{lem555} (ii), we have the following formal
identity
$$\bar{\Phi}'_{\theta^{\alpha},\mu}(f(z),\bar{f}'(w))=\frac{V_{\alpha}^{\mu}\left(
(\partial^{\beta}\bar{f}(w))
_{|\beta|\leq
|\alpha|},z,w\right)}{D^{2|\alpha|-1}(z,w)},\ {\rm in}\
{\cal M},$$ for any $\alpha \in \N^N$, and
$1\leq \mu \leq c$. Thus, plugging this in
(\ref{eq1973}), we obtain, for
$(z,w)\in {\cal M}$, 
\begin{equation}\label{eq1974} \sum_{i=0}^{k_0} s_i
\left(\left(
\frac{V_{\alpha}^{\mu}\left((\partial^{\beta}\bar{f}(w))_ {|\beta|\leq |\alpha|},z,w\right)}{D^{2|\alpha|-1}(z,w)}
\right)_{1\leq \mu \leq c,|\alpha|\leq l_0},\bar{f}'(w)
\right)\left((\chi\circ f)(z)\right)^i \equiv 0. 
\end{equation} We claim that for $(z,w)\in {\cal M}$

\begin{equation}\label{claim}
s_{k_0}\left(\left( \displaystyle
\frac{V_{\alpha}^{\mu}\left((\partial^{\beta}\bar{f}(w))_
 {|\beta|\leq
|\alpha|},z,w\right)}{D^{2|\alpha|-1}(z,w)}\right)_{1\leq
\mu \leq c,|\alpha|\leq l_0},\bar{f}'(w)\right)\not
\equiv 0.
\end{equation}
Indeed, we have first to notice that, by definition,
the left hand side of (\ref{claim}) is equal to
$$s_{k_0}
\left((\bar{\Phi}'_{\theta^{\alpha},\mu}
(f(z),\bar{f}'(w)))_{1\leq
\mu
\leq c,|\alpha|\leq l_0},
\bar{f}'(w)\right).$$ Denote ${\cal
Q}(\omega,\theta)=s_{k_0}\left(
(\bar{\Phi}'_{\theta^{\alpha}}(\omega,\theta))_{\mu
\leq c,|\alpha|\leq l_0},\theta \right)$.  Assuming
(\ref{claim}) false, we would get ${\cal
Q}(f(z),\bar{f}'(w))\equiv 0$, for $(z,w)\in {\cal
M}$. Since $f$ is nondegenerate, one can easily show
that the rank of the formal holomorphic map ${\cal
M}\ni (z,w)\mapsto (f(z),\bar{f}'(w))\in \C^{2n-c}$ is
$2n-c$. (By the rank of a formal holomorphic mapping
$g(x)=(g_1(x),\ldots,g_k(x))$, we mean its rank in the
quotient field of $\C[[x]]$.) By standard arguments
about formal power series, this implies that
${\cal Q}$ is identically zero as a {\it formal} power
series, and hence, identically zero as a convergent
one.  This contradicts (\ref{eq4.2}) and thus proves
(\ref{claim}). To conclude the proof of Proposition
\ref{propr3}, we observe the following. Since each
$s_i$, $0\leq i\leq k_0$, is a polynomial, one
sees that multiplying (\ref{eq1974}) by enough powers
of
$D(z,w)$, we have
reached the desired conclusion in view (i) and (ii) of
Lemma
\ref{lem555}. This finishes the proof of Proposition
\ref{propr3}.$\Box$
\begin{remark}\label{rk2}  By Proposition
\ref{prop4.1}, Proposition \ref{propr3} also holds for
any function
$\chi$ belonging to the algebraic closure of the field
generated by ${\cal C}$.\end{remark}
\begin{remark}\label{rk2bis}
It is worth mentioning that if, in Proposition
\ref{propr3}, $M'$ is furthermore assumed to be
holomorphically nondegenerate, then one can obtain a
more precise statement. Indeed, when $M'$ is
holomorphically nondegenerate, by Theorem \ref{thmir},
the algebraic closure of the field generated by the
family ${\cal C}$ coincides with ${\cal F}_n$. Thus,
in view of Remark \ref{rk2}, we can
apply Proposition
\ref{propr3} to the algebraic functions
$\chi (\omega)=\omega_i$, for $i=1,\ldots,n$, taken as
coordinates. This gives the following proposition.
(Recall also that by Remark \ref{rk1}, when $M'$ is
holomorphically nondegenerate, $l_0=l(M')$, the
Levi-type of
$M'$.)
\end{remark}
\begin{proposition}\label{propbis}
 Let $M\subset \C^n$ be a real analytic
generic submanifold through the origin and $M'\subset \C^n$ be a real algebraic generic submanifold through 
the
origin with the same CR dimension.  Let
$f:(M,0)\rightarrow (M',0)$ be a formal nondegenerate
CR map and assume that $M'$ is holomorphically
nondegenerate.  Then, for
$j=1,\ldots,n$, there exists a positive integer $k_j$
(depending only on $M'$) and a family of convergent
power series
$\delta_{i,j}=\delta_{i,j}\left((\Lambda_{\gamma})_{|\gamma|\leq
l(M')},z,w\right)\in
\C\{(\Lambda_{\gamma}-\partial^{\gamma}\bar{f}(0))_{|\gamma|\leq
l(M')},z,w\}$, $i=0,\ldots,k_j$, such that the formal
identity
$$\sum_{i=0}^{k_j}\delta_{i,j}\left((\partial^{\gamma}\bar{f}(w))_{|\gamma|\leq
l(M')},z,w\right) \left(f_j(z)\right)^i =0$$
holds for $(z,w)\in {\cal M}$, with
$\delta_{k_j,j}\left((\partial^{\gamma}\bar{f}(w))_{|\gamma|\leq
l(M')},z,w\right)\not \equiv 0$ on ${\cal M}$.
\end{proposition}

\begin{remark}\label{rk13}
In view of the works of Baouendi, Ebenfelt and Rothschild \cite{BER1, BER3, BER4}, Proposition
\ref{propr3} can be viewed as a {\it generalized} reflection identity.  We shall propose in the next section
an algebraic interpretation of this identity, which can be compared to the work of Coupet, Pinchuk and Sukhov
\cite{CPS2}.\end{remark}

\section{A principle of analyticity for formal CR power series}\label{sec5} \setcounter{equation}{0} Throughout
this section, which is independent of \S \ref{sec4},
we shall consider one real analytic generic
submanifold $M$, of CR dimension $N$ and of real
codimension $c$ through the origin in
$\C^n$, $n>1$.  We shall also use the notations introduced for $M$ in \S \ref{sec2}. In particular,
the complexification of $M$ is still denoted by ${\cal
M}$. The purpose of \S \ref{sec5} is to prove the
following principle of convergence for formal CR
functions. Recall that 
by a formal CR
function, we mean a formal holomorphic power
series. 
\begin{theorem}\label{th5.1}
Let $M$ be a real analytic generic submanifold at $0\in \C^n$. Let $h(z)$ be a holomorphic formal power series
in $\C[[z]]$, $z\in \C^n$. Assume that: 
\begin{enumerate}
\item[(i)] $M$ is minimal at 0.
\item [(ii)] there exists a formal power series mapping $X(w)=(X_1(w),\ldots,X_m(w))\in
(\C[[w]])^m$, 
$w\in \C^n$, $X(0)=0$, and a family of convergent power series ${\cal U}_j(X,z,w)\in \C\{X,z,w\}$, $j=0,\ldots,l$,
$l\in \N^*$, such that the relation 
\begin{equation}\label{EQUATION1}
\sum_{j=0}^l{\cal
U}_j(X(w),z,w)\left(h(z)\right)^j=0
\end{equation} holds as a formal
power series identity for $(z,w)\in \cal M$, and such
that
\begin{equation}\label{EQUATION2}
{\cal
U}_l(X(w),z,w)\not
\equiv 0,\ {\rm for}\ (z,w)\in {\cal M}.
\end{equation}
\end{enumerate}
Then $h(z)$ is convergent. 
\end{theorem}

The proof of Theorem \ref{th5.1} will be divided in three distinct steps.  \subsection{Algebraic dependence of
the jets.}\label{ssec5.1} \begin{proposition}\label{prop5.1} Let $M$ be real analytic generic submanifold
through the origin in $\C^n$ and $h(z)$ be a formal holomorphic power series in $z=(z_1,\ldots,z_n)$.  Assume
that $h$ satisfies {\rm (ii)} of Theorem
$\ref{th5.1}$.  Then, for any multi-index $\mu\in
\N^n$, there exists two positive integers $l(\mu)$,
$p(\mu)$, a family of convergent power series ${\cal
U}_{i,\mu}\left((\Lambda_{\gamma})_{|\gamma|\leq
|\mu|},z,w\right)\in
\C\{(\Lambda_{\gamma}-\partial^{\gamma}X(0))_{|\gamma|\leq
|\mu|},z,w\}$, $i=0,\ldots,l(\mu)$, such that the
formal identity $$\sum_{i=0}^{l(\mu)} {\cal
U}_{i,\mu}\left((\partial^{\gamma}X(w))_{|\gamma|\leq
|\mu|},z,w\right)\left(\partial^{\mu}h(z)\right)^i=0$$
holds for
$(z,w)\in {\cal M}$,
 and such that
${\cal
U}_{l(\mu),\mu}\left((\partial^{\gamma}X(w))_{|\gamma|\leq
|\mu|},z,w\right)\not \equiv 0$, for $(z,w)\in {\cal
M}$. Here, $X(w)$ is the formal power series mapping
given by {\rm (ii)} of Theorem
$\ref{th5.1}$.
\end{proposition}
{\it Proof.} Since $h(z)$ satisfies (ii)
of Theorem
\ref{th5.1}, there exists a formal power series
mapping $X(w)=(X_1(w),\ldots,X_m(w))\in (\C[[w]])^m$,
$w\in \C^n$, $X(0)=0$, and a family of
convergent power series ${\cal U}_j(X,z,w)\in
\C\{X,z,w\}$, $j=0,\ldots,l$, such that
the formal identities
(\ref{EQUATION1})  and
(\ref{EQUATION2}) hold. For the
proof of the proposition, we assume that the
complexification of $M$ is given by ${\cal
M}=\{(z,w)\in (\C^{2n},0) :  w^*=\bar{\Phi} (z,w')\}$,
with $\bar\Phi$ as defined in \S \ref{ssec2.2}. 
(Recall that
$w=(w',w^*)\in
\C^N\times \C^c$.)  It will be convenient to introduce
for any integer $k$, a subring ${\cal A}_k\subset
\C[[z,w']]$ which is defined as follows.  Let
$$\Pi_X^k
:\C\{(\Lambda_{\gamma}-\partial^{\gamma}X(0))_{|\gamma|\leq
k},z,w'\}\rightarrow \C[[z,w']]$$ be the substitution
homomorphism defined by
$$(\Lambda_{\gamma})_{|\gamma|\leq k}\mapsto
\left(\partial^{\gamma}X(w',\bar{\Phi}(z,w'))\right)
_{|\gamma|\leq k},\ z\mapsto z,\ w'\mapsto w'.$$
${\cal A}_k$ is, by definition, the ring image
$\Pi_X^k (\C\{(\Lambda_{\gamma}-\partial^{\gamma}X(0))
_{|\gamma|\leq k},z,w'\})$.  Finally, we define ${\cal
B}_k$ to be the quotient field of ${\cal A}_k$.   The
reader can now easily check that, to prove the
proposition, it is equivalent to prove that
\begin{equation}\label{stat}
 {\rm (*)}_{\mu}\quad \forall \mu \in \N^{n},\ 
\partial^{\mu}h(z)\ {\rm is\ algebraic\ over\ the\
field}\ {\cal B}_{|\mu|}.\nonumber
\end{equation} We shall prove ${\rm (*)}_{\mu}$ by
induction on $|\mu|$.  For
$|\mu| =0$,  ${\rm (*)}_{0}$ follows
from (\ref{EQUATION1}) and (\ref{EQUATION2}) and the
definition of ${\cal B}_0$.  Assume that
${\rm (*)}_{\mu}$ holds for all $|\mu|=k$.  This means
precisely that, for any $\mu
\in \N^n$ such that $|\mu|=k$, there exist two
positive integers $l(\mu)$, $p(\mu)$, a family of
convergent power series ${a}
_{i,\mu}\left((\Lambda_{\gamma})_{|\gamma|\leq
|\mu|},z,w'\right)\in
\C\{(\Lambda_{\gamma}-\partial^{\gamma}X(0))_{|\gamma|\leq
|\mu|},z,w'\}$, $i=0,\ldots,l(\mu)$, such that the
formal identity \begin{equation}\label{eqMI}
\sum_{i=0}^{l(\mu)}
a_{i,\mu}\left((\partial^{\gamma}X(w',\bar{\Phi}(z,w')))_{|\gamma|\leq
|\mu|},z,w'\right)\left(\partial^{\mu}h(z)\right)^i\equiv
0
\end{equation} holds in $\C[[z,w']]$, and such that
\begin{equation}\label{eqMIR}
a_{l(\mu),\mu}\left((\partial^{\gamma}X(w',\bar{\Phi}(z,w')))_{|\gamma|\leq
|\mu|},z,w'\right)\not \equiv 0.
\end{equation} Moreover, we can choose $l(\mu)$ minimal satisfying a non-trivial relation such as
(\ref{eqMI}).  This implies that \begin{equation}\label{eqRA} \sum_{j=1}^{l(\mu)}
ja_{j,\mu}\left((\partial^{\gamma}X(w',\bar{\Phi}(z,w')))_{|\gamma|\leq
|\mu|},z,w'\right)\left(\partial^{\mu}h(z)\right)^{j-1}\not
\equiv 0, \end{equation} in $\C[[z,w']]$.  In what
follows, for $j=1,\ldots,n$, $1_j$ is the multiindex of
$\N^n$ with $1$ at the $j$-th digit and 0 elsewhere.  Applying $\partial_{z_j}$
for $j=1,\ldots,n$ to (\ref{eqMI}), we obtain
\begin{equation}\label{eq5.3}
\left(\partial^{\mu+1_j}h(z)\right)\sum_{j=1}^{l}
ja_{j,\mu}\left((\partial^{\gamma}X(w',\bar{\Phi}(z,w')))_{|\gamma|\leq
|\mu|},z,w'\right)\left(\partial^{\mu}h(z)\right)^{j-1}\in
{\cal A}_{k+1}[\partial^{\mu}h(z)], \end{equation}
where ${\cal A}_{k+1}[\partial^{\mu}h(z)]$ is the
subring of
$\C[[z,w']]$ generated by ${\cal A}_{k+1}$ and $\partial^{\mu}h(z)$.  By (\ref{eqRA}) and (\ref{eq5.3}), we
see that $\partial^{\mu+1_j}h(z)$ is algebraic over the field ${\cal
B}_{k+1}(\partial^{\mu}h(z))$, which is the subfield
of Frac $\C[[z,w']]$ generated by ${\cal B}_{k+1}$ and
$\partial^{\mu}h(z)$.  Since ${\rm (*)}_{\mu}$ holds,
$\partial^{\mu}h(z)$ is algebraic over ${\cal
B}_k\subset {\cal B}_{k+1}$, and thus, we see that
$\partial^{\mu+1_j}h(z)$ is algebraic over 
${\cal B}_{k+1}$ according to the transitivity of
algebraicity over fields \cite{ZS}. This shows that
${\rm (*)}_{\nu}$ holds for all multiindeces $\nu \in
\N^n$ such that $|\nu|=k+1$. This completes the proof
of
${\rm (*)}_{\mu}$ for all multiindeces $\mu \in \N^n$
and thus the proof of Proposition
\ref{prop5.1}.$\Box$

\subsection{Non-trivial relations at the level of the Segre sets.}\label{ssec5.2} In this section, we shall
make use of the Segre sets mappings $v_j$, $j\in \N$,
associated to $M$ as defined by (\ref{eqss1}) and
(\ref{eqss2}). We shall also keep the notations
introduced in \S \ref{sec2} for $M$.  Our main
purpose here is to establish the following result.
\begin{proposition}\label{prop5.2} Under the
assumptions and notations of Theorem
$\ref{th5.1}$, the following holds. For any multi-index
$\mu\in \N^n$, and for any $d\in \N$, there exist two
positive integers
$\tau=\tau(\mu,d),\ p=p(\mu,d)$, and a family of
convergent power series
$g_{i\mu d}=g_{i\mu
d}\left((\Lambda_{\gamma})_{|\gamma|\leq
p},z,w\right)\in
\C\{(\Lambda_{\gamma}-\partial^{\gamma}X(0))_{|\gamma|\leq
p},z,w\}$, $i=0,\ldots,\tau$, such that the formal
identity
$$\sum_{j=0}^{\tau}g_{j\mu
d}\left((\partial^{\gamma}X\circ
\bar{v}_d)_{|\gamma|\leq p},v_{d+1},\bar{v}_d \right)
\left(\partial^{\mu}h\circ v_{d+1}\right)^j\equiv 0$$
holds in the ring of formal power series in $(d+1)N$
indeterminates, and such that
$g_{\tau \mu d}\left((\partial^{\gamma}X\circ
\bar{v}_d)_{|\gamma|\leq
p},v_{d+1},\bar{v}_d\right) \not
\equiv 0$. Here,
$X(w)$ is the formal power series mapping given by {\rm (ii)} of Theorem
$\ref{th5.1}$ and $N$ is the CR dimension of $M$. 
\end{proposition} 
\begin{remark}\label{rk3} If $M$ is a
generic real analytic submanifold through the origin in
$\C^n$, and
$h(z)$ is a formal holomorphic power series in $z\in
\C^n$ satisfying (ii) of Theorem \ref{th5.1}, then, by
applying Proposition \ref{prop5.1}, for any
multi-index $\mu\in \N^n$, there exist two positive
integers $l(\mu)$, $p(\mu)$, a family of convergent
power series ${\cal
U}_{i,\mu}\left((\Lambda_{\gamma})_{|\gamma|\leq
|\mu|},z,w\right)\in
\C\{(\Lambda_{\gamma}-\partial^{\gamma}X(0))_{|\gamma|\leq
|\mu|},z,w\}$, $i=0,\ldots,l(\mu)$, such that the
formal identity
\begin{equation}\label{eqrem}
\sum_{i=0}^{l(\mu)} {\cal
U}_{i,\mu}\left((\partial^{\gamma}X(w))_{|\gamma|\leq
|\mu|},z,w\right)\left(\partial^{\mu}h(z)\right)^i=0
\end{equation} holds for $(z,w)\in {\cal M}$,
 and such that
${\cal
U}_{l(\mu),\mu}\left((\partial^{\gamma}X(w))_{|\gamma|\leq
|\mu|},z,w\right) \not \equiv 0$  on ${\cal M}$. If,
furthermore, $M$ is minimal at 0, then for $d_0$ large
enough, it follows from Theorem \ref{th2.1} and the
definition of the Segre sets mappings given by
(\ref{eqss1}) and (\ref{eqss2}) that, for
$d\geq d_0$, the holomorphic map $$(\C^{dN},0)\ni
(t_1,t_2,\ldots,t_{d})
\mapsto (v_{d+1}(t_1,t_2,\ldots),
\bar{v}_{d}(t_2,\ldots))\in {\cal M}$$
is generically submersive. Thus, by
elementary facts about formal power series, this
implies that for $d\geq d_0$,
$${\cal U}_{l(\mu),\mu}\left((\partial^{\gamma}X\circ
\bar{v}_d)_{|\gamma|\leq
|\mu|},v_{d+1},\bar{v}_d\right)\not
\equiv 0.$$ This means that the algebraic relations
$$\sum_{i=0}^{l(\mu)} {\cal
U}_{i,\mu}\left((\partial^{\gamma}X\circ
\bar{v}_d)_{|\gamma|\leq
|\mu|},v_{d+1},\bar{v}_d\right)\left(\partial^{\mu}h\circ
v_{d}\right)^i=0
$$ 
will still be non-trivial for $d\geq d_0$. This proves
Proposition
\ref{prop5.2} for $d\geq d_0$. However, in general, 
plugging $z=v_{d+1}$ and $w=\bar{v}_d$ in
(\ref{eqrem}) for $d<d_0$, could lead to trivial
relations. Thus, one has to work a little bit more to
prove Proposition \ref{prop5.2} for $d<d_0$.
\end{remark}

For the proof of Proposition \ref{prop5.2}, we need to introduce the following definition, in which only
the generic submanifold $M$ is involved.
\begin{definition}\label{def11}
Let $M$ be a generic real analytic submanifold through the
origin, of CR dimension $N$, and $v_k$, $k\in \N$, the
associated Segre sets mappings as defined by
(\ref{eqss1}) and (\ref{eqss2}). Let
$Y(w)=(Y_1(w),\ldots,Y_r(w))\in (\C[[w]])^r$, be a
formal power series mapping in
$w=(w_1,\ldots,w_n)$. Given $d\in \N$ and a formal
power series
$q(z)\in \C[[z_1,\ldots,z_n]]$, we say that $q$
{\it satisfies property ${\cal P}(M,Y,d)$} if there
exists a family of convergent power series
$A_j(\Lambda_0,z,w)\in
\C\{(\Lambda_0-Y(0),z,w\}$,
$j=0,\ldots,p$, $p\in \N^*$, such that the identity
$$\sum_{j=0}^pA_j(Y\circ \bar{v}_d,v_{d+1},\bar{v}_d)
\left(q\circ v_{d+1}\right)^j\equiv 0$$
holds in the ring of formal power series in $(d+1)N$
indeterminates and such that $A_p(Y\circ
\bar{v}_d,v_{d+1},\bar{v}_d)\not \equiv 0$.
\end{definition}
We will need the following lemma to derive Proposition \ref{prop5.2}. 

\begin{lemma}\label{lem5.2} Let $M$ be a real analytic generic submanifold through the origin in $\C^n$.  Let
$Y(w)=(Y_1(w),\ldots,Y_r(w))\in (\C[[w]])^r$, be a formal power series mapping in
$w=(w_1,\ldots,w_n)$.  Let $d\in \N$ and $q(z)\in 
\C[[z]]$, $z=(z_1,\ldots,z_n)$.  Then, if $q(z)$
satisfies property ${\cal P}(M,Y,d+2)$, there exists an
integer
$n_0$ (depending on $Y$, $q$ and $d$) such that $q(z)$
satisfies property ${\cal
P}(M,(\partial^{\beta}Y)_{|\beta|\leq n_0},d)$. 
\end{lemma}
{\it Proof of Lemma $\ref{lem5.2}$.} Let
$Y$,
$q(z)$ and $d$ be as in the Lemma.  We assume that
$q(z)$ satisfies property ${\cal P}(M,Y,d+2)$. By
definition, there exists a family of convergent power
series
$A_j(\Lambda_0,z,w)\in \C\{(\Lambda_0-Y(0),z,w\}$,
$j=0,\ldots,p$, $p\in \N^*$, such that the formal
identity
\begin{equation}\label{eqRR1} \sum_{j=0}^p
A_j(Y\circ \bar{v}_{d+2},v_{d+3},\bar{v}_{d+2})
\left(q\circ
v_{d+3}\right)^j\equiv 0
\end{equation} holds and such that $A_p(Y\circ
\bar{v}_{d+2},v_{d+3},\bar{v}_{d+2})\not
\equiv 0$. Here, 
$$v_{d+3}=v_{d+3}(t_1,t_2,\ldots,t_{d+3})=v_{d+3}(t_1,t')=v_{d+3}(t),$$
$$\bar{v}_{d+2}=\bar{v}_{d+2}(t_2,\ldots,t_{d+3})=\bar{v}_{d+2}(t').$$ Thus, (\ref{eqRR1}) holds in the ring
$\C[[t_1,\ldots,t_{d+3}]]$. For simplicity of notations, we put, for $j=0,\ldots,p$,
\begin{equation}\label{eqr9}
\Theta_j (t)=A_j\left(((Y\circ \bar{v}_{d+2})(t'),
v_{d+3}(t),\bar{v}_{d+2}(t')\right).
\end{equation}
Thus, (\ref{eqRR1}) can be rewritten as
\begin{equation}\label{eqr10}
\sum_{j=0}^{p}\Theta_j (t)\left((q \circ
v_{d+3})(t)\right)^j\equiv 0,\ {\rm with}
\end{equation}

\begin{equation}\label{eqr11} \Theta_{p}(t)\not \equiv 0,\ {\rm in}\ \C[[t]].  \end{equation} Consider the set
${\cal E}$ defined by $$\{\alpha \in \N^N:\exists j\in
\{1,\ldots,p\},\ {\rm such\ that}\
\left[\frac{\partial^{|\alpha|} \Theta_j}{\partial t_1^{\alpha}}(t)\right]_{t_1=t_3}\not \equiv 0,\ {\rm in}\
\C[[t']]\}.$$ Observe that by (\ref{eqr11}), there exists a multiindex $\alpha \in \N^N$ such that
$$\left[\frac{\partial^{|\alpha|} \Theta_p}{\partial t_1^{\alpha}}(t)\right]_{t_1=t_3}\not \equiv 0$$ in
$\C[[t']]$.  This implies that ${\cal E}$ is not empty.  Let $\alpha^0 \in \N^N$ such that $|\alpha^0|={\rm
min}\{|\beta|:\beta \in {\cal E}\}$.  Then, if we
apply $\displaystyle \frac{\partial^
{|\alpha^0|}}{\partial t_1^{\alpha^0}}$ to
(\ref{eqr10}), it follows  from Leibniz's formula that
\begin{equation}\label{eq5.6}
\frac{\partial^{|\alpha^0|}\Theta_0}{\partial t_1^{\alpha^0}}(t)+ \sum_{j=1}^{p}
\frac{\partial^{|\alpha^0|}\Theta_j}{\partial
t_1^{\alpha^0}}(t)\left((q\circ v_{d+3})(t)\right)^j=
\sum_{\beta
\in \N^N, |\beta|<|\alpha^0| \atop 1\leq j\leq p}
\frac{\partial^{\beta}\Theta_j}{\partial t_1
^{\beta}}(t)\vartheta_{\beta,j}(t), \end{equation}
where, for any
$\beta$, $j$, $\vartheta_{\beta,j}(t) \in \C[[t]]$. By
the choice of $\alpha^0$, we have, for
$|\beta|<|\alpha^0|$, \begin{equation}\label{eqr12}
\left[\frac{\partial^{\beta}\Theta_j}{\partial t_1^{\beta}}(t_1,t')\right]_{t_1=t_3}\equiv 0, \ {\rm in}\
\C[[t']],\ j=1,\ldots,p.  \end{equation} 
Thus, if we restrict equation (\ref{eq5.6}) to $t_1=t_3$, we get by
(\ref{eqr12}) the following identity in the ring $\C[[t']]$ 
\begin{equation}\label{eqr14} \sum_{j=0}^{p}
\left[\frac{\partial^{|\alpha^0|}\Theta_j}{\partial t_1^{\alpha^0}}(t_1,t')\right]_{t_1=t_3}
\left((q \circ
v_{d+3})(t_3,t_2,t_3,\ldots,t_{d+3})\right)^j\equiv 0. 
\end{equation} Here again, for simplicity of
notations, we put \begin{equation}\label{eq2001}
\Theta'_j(t')=\left[\frac{\partial^{|\alpha^0|}\Theta_j}{\partial
t_1^{'\alpha^0}}(t_1,t')\right]_{t_1=t_3}. 
\end{equation}  We observe that, by the choice of
$\alpha^0$, there exists $j\in \{1,\ldots,p\}$ such
that $\Theta'_j(t')\not \equiv 0$.  Denote $m_1={\rm
Sup}\{j\in
\{1,\ldots,p\}:\Theta'_j(t')\not \equiv 0\}$.  It
follows from the reality  condition (\ref{eq2000}) and
the definition of the Segre sets mappings given by
(\ref{eqss1}) and (\ref{eqss2}) that
$$v_{d+3}(t_3,t_2,t_3,\ldots,t_{d+1})=v_{d+1}(t_3,\ldots,t_{d+3}).$$ Thus, (\ref{eqr14}) reads
as \begin{equation}\label{eqr15} \sum_{j=0}^{m_1}
\Theta'_j(t')\left((q \circ
v_{d+1})(t_3,\ldots,t_{d+3})\right)^j\equiv 0,
\end{equation} with, moreover,
\begin{equation}\label{ND}
\Theta'_{m_1}(t')\not \equiv 0.
\end{equation}

 {\it First case
:}  $d\geq 1$.  (\ref{ND}) implies that there
exists $\beta^0 \in \N^N$ such that
\begin{equation}\label{eqb1}
\left[\frac{\partial^{|\beta^0|}\Theta'_{m_1}}{\partial
t_2^{\beta^0}}(t_2,t_3,\ldots,t_{d+3})\right]_{t_2=t_4}\not \equiv 0.  \end{equation} Thus, applying
$\displaystyle \frac{\partial^{|\beta^0|}}{\partial t_2^{\beta^0}}$ to (\ref{eqr15}) and after evaluation at
$t_2=t_4$, we obtain in the ring
$\C[[t_3,t_4,\ldots,t_{d+3}]]$
\begin{equation}\label{eqr16}
\sum_{j=0}^{m_1}\left[\frac{\partial^{|\beta^0|}\Theta'_j}{\partial
t_2^{\beta^0}}(t')\right]_{t_2=t_4}\left((q\circ
v_{d+1})(t_3,\ldots,t_{d+3})\right)^j\equiv 0.
\end{equation} We shall now see that (\ref{eqr16}) gives the statement of the Lemma.  By definition of the
Segre sets mappings given by (\ref{eqss1}) and
(\ref{eqss2}) and by the definition of the
$\Theta_j$ given by (\ref{eqr9}), we have, for
$0\leq j\leq m_1$, \begin{eqnarray}
\Theta_j(t_1,t')&=& A_j\left(Y\circ
\bar{v}_{d+2}(t'),t_1,\Phi
(\bar{v}_{d+2}(t'),t_1),\bar{v}_{d+2}(t')\right)\nonumber
\\ &=&G^1_j\left((Y\circ
\bar{v}_{d+2})(t'),\bar{v}_{d+2}(t'),t_1\right)\nonumber,
\end{eqnarray} where $G^1_j
\left(\Lambda_0,w,t_1\right)\in
\C\{\Lambda_{0}-Y(0),w,t_1\}$.  Using (\ref{eq2001}),
we obtain \begin{eqnarray}
\Theta'_j(t')&=&\left[\frac{\partial^{|\alpha^0|}G^1_j}{\partial
t_1^{\alpha^0}}\right]\left((Y\circ
\bar{v}_{d+2})(t'),\bar{v}_{d+2}(t'),t_3\right)
\nonumber
\\
 &=&\left[\frac{\partial^{|\alpha^0|}G^1_j}{\partial
t_1^{\alpha^0}}\right]\left((Y(t_2,\bar{\Phi}
(v_{d+1},t_2)),t_2,\bar{\Phi}(v_{d+1},t_2),t_3\right).
\nonumber
\end{eqnarray} Here,
$v_{d+1}=v_{d+1}(t_3,t_4,\ldots,t_{d+3})$.  As a consequence, we have
$$\frac{\partial^{|\beta^0|}\Theta'_j}{\partial
t_2^{\beta^0}}(t')=G^2_j\left((\partial^{\gamma}Y(t_2,\bar{\Phi}
(v_{d+1},t_2))_{|\gamma|\leq
|\beta^0|},v_{d+1},t_2,t_3\right),$$ where
$G^2_j=G^2_j\left((\Lambda_{\gamma})_{|\gamma|\leq
|\beta^0|},w,t_2,t_3\right)\in
\C\{(\Lambda_{\gamma}-\partial^{\gamma}Y(0))_{|\gamma|\leq
|\beta^0|},w,t_2,t_3\}$. Here again, by
(\ref{eq2000}), (\ref{eqss1}) and (\ref{eqss2}), we
have
$\bar{v}_{d+2}(t_4,t_3,t_4,t_5,\ldots)
=\bar{v}_d(t_4,t_5,\ldots)$,
 and thus \begin{eqnarray}
\left[\frac{\partial^{|\beta^0|}\Theta'_j}{\partial
t_2^{\beta^0}}(t')\right]_{t_2=t_4}&=&
G^2_j\left(((\partial^{\gamma}Y\circ \bar{v}_{d+2})
(t_4,t_3,t_4,\ldots))_{|\gamma|\leq
|\beta^0|},v_{d+1}(t_3,t_4,\ldots),t_4,t_3\right)\nonumber
\\ &=&G^2_j\left(((\partial^{\gamma}Y\circ
\bar{v}_{d})(t_4,\ldots,t_{d+3}))_{|\gamma|\leq
|\beta^0|},v_{d+1}(t_3,t_4,\ldots),t_4,t_3\right)\nonumber
\\ &=& B_j\left((\partial^{\gamma}Y\circ
\bar{v}_d)_{|\gamma|\leq
|\beta^0|},v_{d+1},\bar{v}_d\right),  \label{eqf1} 
\end{eqnarray} 
where $B_j$ for $j=0,\ldots,m_1$, is a convergent power
series in its arguments. Consequently, from
(\ref{eqr16}) and (\ref{eqf1}) we have the relation
$$\sum_{j=0}^{m_1}B_j\left((\partial^{\gamma}Y\circ
\bar{v}_d)_{|\gamma|\leq
|\beta^0|},v_{d+1},\bar{v}_d\right) \left(q\circ
v_{d+1}\right)^j\equiv 0,\ {\rm in}\
\C[[t_3,\ldots,t_{d+3}]],$$ which is non-trivial
according to (\ref{eqb1}) and (\ref{eqf1}). In
conclusion,
$q(z)$ satisfies property ${\cal
P}(M,(\partial^{\gamma}Y)_{|\gamma|\leq |\beta^0|},d)$.

{\it Second case :} $d=0$. In this case, almost the same procedure used in the case $d\geq 1$ can be applied.
 Indeed, by (\ref{ND}), we have
$\Theta'_{m_1}(t_2,t_3)\not
\equiv 0$ and therefore, there exists a multi-index
$\varrho^0\in
\N^N$ such that

\begin{equation}\label{eqb2} \left[\frac{\partial^{|\varrho^0|}\Theta'_{m_1}}{\partial
t_2^{\varrho^0}}(t_2,t_3)\right]_{t_2=0}\not \equiv 0.  \end{equation} Thus, applying $\displaystyle
\frac{\partial^{|\varrho^0|}}{\partial t_2^{\varrho^0}}$ to (\ref{eqr15}) and after evaluation at $t_2=0$, we
obtain in the ring $\C[[t_3]]$ \begin{equation}\label{eqr17}
\sum_{j=0}^{m_1}\left[\frac{\partial^{|\varrho^0|}\Theta'_j}{\partial t_2^{\varrho^0}}(t_2,t_3)\right]_{t_2=0}
\left((q\circ v_{1})(t_3)\right)^j\equiv 0. 
\end{equation}  As in the case $d\geq 1$, we have for
$j=1,\ldots,m_1$,
$$\frac{\partial^{|\varrho^0|}\Theta'_j}{\partial t_2^{\varrho^0}}(t_2,t_3)=
G^3_j\left((\partial^{\gamma}Y (t_2,\bar{\Phi}
(v_{1}(t_3),t_2))_{|\gamma|\leq
|\varrho^0|},v_{1}(t_3),t_2,t_3\right),$$ where
$G_j^3=G^3_j\left((\Lambda_{\gamma})_{|\gamma|\leq
|\varrho^0|},w,t_2,t_3\right)\in
\C\{(\Lambda_{\gamma}-\partial^{\gamma}Y
(0))_{|\gamma|\leq |\varrho^0|},w,t_2,t_3\}$.  By the
normality of the coordinates for $M$, we have for
$j=0,\ldots,m_1$, 
\begin{equation}\label{eq2002} 
\left[\frac{\partial^{|\varrho^0|}\Theta'_j}{\partial
t_2^{\varrho^0}}(t_2,t_3)\right]_{t_2=0}=
G^3_j\left((\partial^{\gamma}Y (0))_{|\gamma|\leq
|\varrho^0|},v_{1}(t_3),0,t_3\right).  \end{equation}
We leave it to the reader to check that, similarly to
the case
$d\geq 1$, (\ref{eq2002}), (\ref{eqr17}) and (\ref{eqb2}) give the desired statement of the lemma for
$d=0$, i.e. that $q(z)$ satisfies property ${\cal
P}(M,(\partial^{\gamma}Y)_{|\gamma|\leq |\beta^0|},0)$.
This completes the proof of Lemma
\ref{lem5.2}.$\Box$

{\it Proof of Proposition $\ref{prop5.2}$.}
Let $\mu \in \N^n$.  Since $h(z)$ satisfies (ii) of
Theorem \ref{th5.1}, by Proposition \ref{prop5.1},
there exist two positive integers $l(\mu)$, $p(\mu)$ 
and a family of  convergent power series ${\cal
U}_{i,\mu}\left((\Lambda_{\gamma})_{|\gamma|\leq
|\mu|},z,w\right)\in
\C\{(\Lambda_{\gamma}-\partial^{\gamma}X(0))_{|\gamma|\leq
|\mu|},z,w\}$, $i=0,\ldots,l(\mu)$, such that the
formal identity
\begin{equation}\label{EQUATION3}
\sum_{i=0}^{l(\mu)}
{\cal
U}_{i,\mu}\left((\partial^{\gamma}X(w))_{|\gamma|\leq
|\mu|},z,w\right)\left(\partial^{\mu}h(z)\right)^i=0
\end{equation} 
holds on
${\cal M}$ and such that 
\begin{equation}\label{EQUATION4}
{\cal
U}_{l(\mu),\mu}\left((\partial^{\gamma}X(w))_{|\gamma|\leq
|\mu|},z,w\right) \not =0
\end{equation} for $(z,w)\in {\cal M}$.
Notice that to prove the proposition we have to show
that for any
$d\in \N$, there exists $p=p(\mu,d)$ such that
$\partial^{\mu}h(z)$ satisfies property ${\cal
P}(M,(\partial^{\beta}X)_{|\beta|\leq p},d)$.  Since
$M$ is minimal at $0\in M$, it follows from  Theorem
\ref{th2.1} and the definition of the Segre sets
mappings given by (\ref{eqss1}) and (\ref{eqss2}) that
there exists
$d_0\in
\N$ (which can be assumed to be even) such that the
holomorphic map 
$ \C^{(d_0+1)N}\ni (t_1,\ldots,t_{d_0+1})\mapsto
(v_{d_0+1}(t_1,t_2,\ldots,t_{d_0+1}),
\bar{v}_{d_0}(t_2,\ldots,t_{d_0+1}))\in
{\cal M}$ is generically submersive.  
By elementary facts about formal  power series and by
(\ref{EQUATION4}), this implies that
$${\cal
U}_{l(\mu),\mu}\left((\partial^{\gamma}X\circ
\bar{v}_{d_0})_{|\gamma|\leq
|\mu|},v_{d_0+1},\bar{v}_{d_0}\right)\not
\equiv 0.$$
This means that the following algebraic relation 
obtained from (\ref{EQUATION3}),
\begin{equation}\label{eqr20}
\sum_{i=0}^{l(\mu)} {\cal
U}_{i,\mu}\left((\partial^{\gamma}X\circ
\bar{v}_{d_0})_{|\gamma|\leq
|\mu|},v_{d_0+1},\bar{v}_{d_0}\right)\left((\partial^{\mu}h)\circ
v_{d_0+1}\right)^i\equiv 0, \end{equation} is
still non-trivial, i.e. that
$\partial^{\mu}h(z)$ satisfies property
${\cal P}(M,(\partial^{\beta}X)_{|\beta|\leq
|\mu|},d_0)$. (Observe that in Remark \ref{rk3}, we
have shown that
$\partial^{\mu}h(z)$ satisfies property ${\cal
P}(M,(\partial^{\beta}X)_{|\beta|\leq |\mu|},d)$ for
any $d\geq d_0$.)  Applying Lemma
\ref{lem5.2} to $q(z)=\partial^{\mu}h(z)$ and
$Y=(\partial^{\beta}X)_{|\beta|\leq |\mu|}$, we obtain
that there exists
$p(\mu,d_0-2)\in \N$ such that $\partial^{\mu}h(z)$
satisfies property ${\cal
P}(M,(\partial^{\beta}X)_{|\beta|\leq
p(\mu,d_0-2)},d_0-2)$. Hence,  using inductively Lemma
\ref{lem5.2}, we obtain that for any even number
$0\leq d\leq d_0$, there exists
$p(\mu,d)\in \N$ such that 
$\partial^{\mu}h(z)$ satisfies property ${\cal
P}(M,(\partial^{\beta}X)_{|\beta|\leq p(\mu,d)},d)$.
Since by Remark \ref{rk3}, 
$\partial^{\mu}h(z)$ satisfies also property ${\cal
P}(M,(\partial^{\beta}X)_{|\beta|\leq |\mu|},d_0+1)$,
we can  again, in the same way, use Lemma
\ref{lem5.2} to conclude that for any odd number
$1\leq d\leq d_0$, there exists
$p(\mu,d)\in \N$ such that $\partial^{\mu}h(z)$
satisfies property ${\cal
P}(M,(\partial^{\beta}X)_{|\beta|\leq  p(\mu,d)},d)$.
This completes the proof of Proposition
\ref{prop5.2}.$\Box$

\subsection{Propagation procedure.}\label{ssec5.3}
We prove here the last proposition needed for the proof
of Theorem \ref{th5.1}.

\begin{proposition}\label{prop5.3}
Let $M$ be a generic real analytic submanifold through the origin, and $v_k$, $k\in \N$, the associated
Segre sets mappings as defined by {\rm (\ref{eqss1})} and {\rm (\ref{eqss2})}. Let $h(z)\in \C[[z_1,\ldots,z_n]]$ and
$d\in \N$. Let $\mu \in \N^n$ and assume that there exists $Y_{\mu
d}(w)$, a formal power
series mapping in $w=(w_1,\ldots,w_n)$, such that
$\partial^{\mu}h(z)$
 satisfies property ${\cal P}(M,Y_{\mu d},d+2)$ as 
defined in Definition $\ref{def11}$. Then, the
following
 holds. If for any multiindices $\nu \in \N^n$,
$\partial^{\nu}h\circ v_{d+1}$ is convergent, then 
 $\partial^{\mu}h\circ v_{d+3}$ is convergent.\end{proposition}

For the proof of Proposition \ref{prop5.3}, 
we need the following two lemmas which are both
consequences of the Artin approximation theorem
\cite{A}. We refer the reader to \cite{M3} for the
proof of the first one and to \cite{BER4, M3} for the
proof of the second one. 
\begin{lemma}\label{lem5.3} 
 Let ${\cal
T}(x,u)=({\cal T}_1(x,u),\ldots,{\cal T}_r(x,u))\in
(\C[[x,u]])^r$,
$x\in \C^q$, $u\in \C^s$, with ${\cal T}(0)=0$.  Assume that
${\cal T}(x,u)$ satisfies an identity in the  ring
$\C[[x,u,y]]$, $y\in \C^q$, of the form
$$\varphi ({\cal T}(x,u);x,u,y)=0,$$
where $\varphi \in
\C[[W,x,u,y]]$ with $W\in \C^r$.  Assume, furthermore,
that for any multi-index $\beta \in \N^q$, the formal power
series $\displaystyle
\left[\frac{\partial^{|\beta|}\varphi}{\partial y^{\beta}}
(W;x,u,y)\right]_{y=x}$ is convergent,
i.e.  belongs to $\C\{W,x,u\}$.  Then, for any
given positive integer $e$, there exists an $r$-tuple of
convergent power series ${\cal T}^e(x,u)\in
(\C\{x,u\})^r$ such that $\varphi ({\cal
T}^e(x,u);x,u,y)=0$ in
$\C[[x,u,y]]$ and such that ${\cal T}^e(x,u)$
agrees up to order $e$ (at 0) with ${\cal T}(x,u)$.
\end{lemma} 

\begin{lemma}\label{lem5.4}
Any formal power series in $r$
indeterminates, which
is algebraic over the field of meromorphic
functions (in $r$ variables), must be 
convergent.
\end{lemma} 
 \begin{remark}\label{rk5} We would like to mention
that the use of the Artin approximation theorem is not
a novelty in the study of many mapping problems ({\it
cf.}  \cite{D, BR, M2, BER4, M3} as well as many other
articles). \end{remark}
{\it Proof of Proposition $\ref{prop5.3}$.}
Let $\mu \in \N^n$ and $d\in \N$ be as in the statement
of the proposition. We assume that for any multiindex
$\nu \in \N^n$,
$\partial^{\nu}h\circ v_{d+1}$ is convergent. By assumption, there exists $Y_{\mu
d}(w)\in
(\C[[w]])^r$, $r=r(\mu,d)\in \N^*$, a formal power
series mapping in $w=(w_1,\ldots,w_n)$ such that
$\partial^{\mu}h(z)$
 satisfies property ${\cal P}(M,Y_{\mu d},d+2)$. By
definition, this means that there exists a
 family of convergent power series
$A_j=A_j(\Lambda_0,z,w)\in \C\{\Lambda_0-Y_{\mu
d}(0),z,w\}$, $j=0,\ldots,k$, 
 $k=k(\mu,d)$, such that the formal power series identity
 \begin{equation}\label{eqfinal}
 \sum_{j=0}^kA_j\left((Y_{\mu d}\circ
\bar{v}_{d+2})(t'),v_{d+3}(t),\bar{v}_{d+2}(t')\right)
\left(((\partial^{\mu}h)\circ
v_{d+3})(t)\right)^j\equiv 0
 \end{equation}
 holds in $\C[[t]]$, with
$t=(t_1,t')=(t_1,t_2,\ldots,t_{d+3})$ and such that 
 \begin{equation}\label{eqfinale}
 A_{k}((Y_{\mu d}\circ
\bar{v}_{d+2})(t'),v_{d+3}(t),\bar{v}_{d+2}(t'))\not
\equiv 0.
\end{equation}
We would like to apply Lemma \ref{lem5.3} with
$y=t_1$, $x =t_3$, 
$u=(t_2,t_4,t_5,\ldots,t_{d+3})$, 
${\cal T}(x,u)=(Y_{\mu d}\circ
\bar{v}_{d+2})(t_2,\ldots,t_{d+3})-Y_{\mu d}(0)$, 
$W=\Lambda'_0$ ($\Lambda'_0 \in \C^r$) and

\begin{multline}\label{eq2004}
\varphi
(\Lambda'_0;t_3,(t_2,t_4,t_5,\ldots,t_{d+3}),t_1)=\\
\sum_{j=0}^kA_j\left(\Lambda'_0+Y_{\mu
d}(0),v_{d+3}(t),\bar{v}_{d+2}(t')\right)
\left(((\partial^{\mu}h)\circ v_{d+3})(t)\right)^j.
\end{multline}
For this, one has to check that any derivative with respect to $t_1$ of
$\varphi$ evaluated at $t_1=t_3$ is in fact convergent
with respect to the variables 
$\Lambda'_0$ and $t'$. Because of the analyticity of
the functions $A_i$, $i=0,\ldots,k$  (and of the
Segre sets mappings), we see that we have only to
consider the derivatives of 
$[\partial^{\mu}h\circ v_{d+3}(t)]^j$,
for $j=0,\ldots,l$, evaluated at $t_1=t_3$. These
derivatives involve  analytic terms coming for the
differentiation of
$v_{d+3}$ (which are convergent) and products involving
powers of derivatives of
$h$ evaluated at
$t_1=t_3$. Let $[(\partial^{\gamma}h)\circ
v_{d+3}(t)]_{t_1=t_3}$ be such a derivative for some
$\gamma \in \N^n$. By the reality
condition (\ref{eq2000}) and by (\ref{eqss1}) and
(\ref{eqss2}), we have
$$v_{d+3}(t_3,t_2,t_3,t_4,\ldots,t_{d+3})=v_{d+1}(t_3,t_4,\ldots,t_{d+3}).$$
Thus, $[(\partial^{\gamma}h)\circ
v_{d+3}(t)]_{t_1=t_3}=\left((\partial^{\gamma}h)\circ
v_{d+1}\right)(t_3,\ldots,t_{d+3})$ which is convergent
by our hypothesis. As a consequence,
$\varphi$ satifies the assumptions of Lemma
\ref{lem5.3}. Thus, by applying that Lemma, one obtains
for any positive integer $e$,  a convergent power
series mapping ${\cal T}^e(t')$, which agrees up to
order $e$ with 
$(Y_{\mu d} \circ \bar{v}_{d+2})(t')$ and such that
$$\sum_{j=0}^k A_j\left({\cal
T}^e(t'),v_{d+3}(t),\bar{v}_{d+2}(t')\right) 
\left(((\partial^{\mu}h)\circ
v_{d+3})(t)\right)^j\equiv 0.$$ Observe that
(\ref{eqfinale}) implies that, for
$e$ large enough, say $e=e_0$, the following condition
will be satisfied 
$$A_k\left({\cal
T}^{e_0}(t'),v_{d+3}(t),\bar{v}_{d+2}(t')\right)\not
\equiv 0,$$ in $\C[[t]]$. This allows one to apply
Lemma
\ref{lem5.4} to conclude that 
$\partial^{\mu}h\circ v_{d+3}$ is
convergent.$\Box$

\subsection{Completion of the proof of Theorem
\ref{th5.1}}
Let
$h(z)$ be the formal power series of the Theorem and
$X(w)$ the associated formal power series mapping
given by (ii).  By Proposition
\ref{prop5.2}, for any multi-index $\mu\in \N^n$, and for any $d\in \N$, there exists two positive integers
$\tau=\tau(\mu,d),\ p=p(\mu,d)$, and a family of convergent power series
$g_{i\mu d}=g_{i\mu
d}\left((\Lambda_{\gamma})_{|\gamma|\leq
p},z,w\right)\in
\C\{(\Lambda_{\gamma}-\partial^{\gamma}X(0))_{|\gamma|\leq
p},z,w\}$, $i=0,\ldots,\tau$, such that the identity
$$\sum_{j=0}^{\tau}g_{j\mu
d}\left((\partial^{\gamma}X\circ
\bar{v}_d)_{|\gamma|\leq p},v_{d+1},\bar{v}_d\right)
\left((\partial^{\mu}h)\circ v_{d+1}\right)^j\equiv 0$$
holds in $\C[[t]]$ where $t=(t_1,\ldots,t_{d+1})\in
\C^{(d+1)N}$, and with the additional nondegeneracy
condition $g_{\tau \mu d}\left((\partial^{\gamma}X\circ
\bar{v}_d)_{|\gamma|\leq
p},v_{d+1},\bar{v}_d\right)\not
\equiv 0$. In view of Definition \ref{def11}, this means that, for any multiindex $\mu \in
\N^n$ and for any $d\in \N$,
$\partial^{\mu}h(z)$ satisfies property
${\cal P}(M,(\partial^{\gamma}X)_{|\gamma|\leq p(\mu,d)},d)$. Observe first that 
since for any $\nu \in \N^n$, $\partial^{\nu}h(z)$
satisfie property
${\cal P}(M,(\partial^{\gamma} X)_{|\gamma|\leq
p(\nu,0)},0)$, it  follows from Lemma \ref{lem5.4}
 (and from Definition \ref{def11}) that
$\partial^{\nu}h\circ v_1$ is convergent for any
multiindex
$\nu \in \N^n$. From this and the fact that 
$\partial^{\mu}h(z)$ satisfies property ${\cal
P}(M,(\partial^{\gamma}X)_{|\gamma|\leq p(\mu,2)},2)$,
it follows from Proposition
\ref{prop5.3} that $\partial^{\mu}h\circ v_{3}$ is convergent, for all multiindices $\mu \in \N^n$. Thus, by induction,
we see that Proposition \ref{prop5.3} gives that for any odd number $d$, and for any multiindex 
$\mu \in \N^n$,
$\partial^{\mu}h\circ v_d$ is convergent. Choose
$d_0\in
\N$ satisfying the statement of Theorem
\ref{th2.1}. Without loss of generality,
$d_0$ can be assumed to be odd. By the previous considerations, we know that $h\circ v_{d_0}$ is convergent in
some neighborhood $U$ of $0\in \C^{d_0N}$. By 
Theorem \ref{th2.1}, there exists
$(t_1^0,\ldots,t_{d_0}^0)\in U$ such that $v_{d_0}(t_1^{0},\ldots,t_{d_0}^0)=0$ and such that
$v_{d_0}$ is submersive at $(t_1^0,\ldots,t_{d_0}^0)$. Thus, we may apply the rank theorem to conclude that
$v_{d_0}$ has a right convergent inverse $\theta (z)
\in (\C\{z\})^{d_0N}$ defined near $0\in \C^n$ such
that
$\theta (0)=(t_1^0,\ldots,t_{d_0}^0)$  and such that $v_{d_0}\circ \theta (z)=z$. This implies that $h(z)$
is convergent. The proof of Theorem \ref{th5.1} is
complete.

\section{Proofs of Theorem \ref{th1.1} and  Theorem \ref{cor1.2}}\label{sec6}
\setcounter{equation}{0} 
{\it Proof of Theorem $\ref{th1.1}$.} 
Recall that $M'$ is given near 0 by (\ref{eqM"})
and  that
${\cal C}$ is the family of algebraic holomorphic functions (in $n$ variables)
defined by (\ref{eqSVM}) and constructed from $M'$. 
By Proposition \ref{propr3},
for any $\chi \in {\cal C}$, there exists $l_0\in \N^*$, a positive integer $k_0$ and a family of convergent
power series
$\delta_i=\delta_i\left((\Lambda_{\gamma})_{|\gamma|\leq
l_0},z,w\right)\in
\C\{(\Lambda_{\gamma}-\partial^{\gamma}\bar{f}(0))
_{|\gamma|\leq
l_0},z,w\}$, $i=0,\ldots,k_0$, such that the formal
power series identity

\begin{equation}\label{gd}
\sum_{i=0}^{k_0}\delta_i\left((\partial^{\gamma}\bar{f}(w))_{|\gamma|\leq
l_0},z,w\right)\left((\chi\circ f)(z)\right)^i =0,
\end{equation}
holds for
$(z,w)\in {\cal M}$ such that
$\delta_{k_0}\left((\partial^{\gamma}\bar{f}(w))_{|\gamma|\leq
l_0},z,w\right)\not \equiv 0$ in ${\cal M}$.  Since
$M$ is minimal at 0 and $\chi \circ f$ satisfies (ii)
of Theorem \ref{th5.1} by (\ref{gd}), we may apply that
theorem to conclude that $\chi \circ f$ is
convergent.  In other words, for any $\alpha \in \N^N$
and for any
$1\leq
\nu \leq c$,
$q_{\alpha,\nu}\circ f$ is convergent.  To conclude the proof of Theorem \ref{th1.1}, we have to show that
this implies that the reflection mapping
\begin{equation}\label{eqouf} \C^n\times \C^N \ni
(z,\theta)\mapsto
\bar{\Phi}'(f(z),\theta)\in \C^c \end{equation} is 
convergent.  To see this, it suffices to observe that
since $M'$ is real algebraic, the map 
$\C^n \times \C^N \ni
(\zeta,\theta)\mapsto
\bar{\Phi}_{\nu}'(\omega,\theta)\in \C$, 
$1\leq \nu \leq c$, is algebraic, and thus, an
approximation argument similar to the one used in
the proof of Proposition 1 from
\cite{M2} shows that for any $\nu \in \{1,\ldots,c\}$, 
$\bar{\Phi}_{\nu}'(\omega,\theta)$ is algebraic over
the field generated by $\C$ and the family of algebraic
functions ${\cal C}$ and $\theta$.  Since $f$ is
nondegenerate, this implies that the formal power
series $\C^n\times \C^N \ni (z,\theta) \mapsto
\bar{\Phi}_{\nu}'(f(z),\theta)$ is algebraic over the
field generated by $\C$, the family of formal power
series
${\cal C}_f=\left((q_{\beta, \nu}\circ
f)(z)\right)_{\beta
\in
\N^N \atop 1\leq \nu \leq c}$ and $\theta$.  But since
the family ${\cal C}_f$ is a family of convergent
power series, Lemma \ref{lem5.4} implies that the
formal power series $\bar{\Phi}_{\nu}'(f(z),\theta)$
is actually convergent for any $\nu \in
\{1,\ldots,c\}$.  This completes the proof of Theorem
\ref{th1.1}.$\Box$

\begin{remark}\label{rk8} If, in Theorem \ref{th1.1},
the target manifold is given in normal coordinates
i.e. if $\bar \Phi'(\omega,0)=\omega^*$ where $\bar
\Phi'$ is given by (\ref{eqM"}), then the
following holds. The normal components $f^*\in \C^c$
(as in (\ref{split})) of a formal nondegenerate CR map
$f$ from a real analytic generic submanifold into a
real algebraic one are necessarily convergent provided
that the source manifold is minimal. Indeed, this
follows from Theorem \ref{th1.1} by taking
$\theta =0$.
\end{remark}
{\it Proof of Theorem $\ref{cor1.2}$.} By the
Taylor expansion (\ref{taylor}) and by Theorem
\ref{th1.1}, we know that for any $\beta \in \N^N$ and
any $1\leq \nu \leq c$, $q_{\beta, \nu}\circ f$ is
convergent. Equivalently, we have the convergence of
$\chi
\circ f$ for any algebraic function $\chi \in {\cal C}$, 
where ${\cal C}$ is the
family of algebraic functions defined by
(\ref{eqSVM}). Observe that since
$f$ is nondegenerate, Lemma \ref{lem5.4} implies that
for any algebraic holomorphic function $q=q(\omega)$
in the algebraic closure of the field generated by the
family
${\cal C}$, $q\circ f$ must also be convergent. 
To conclude that $f$ is
convergent when $M'$ is holomorphically nondegenerate,
 it suffices to apply Theorem \ref{thmir} (iii) which
states that this algebraic closure, in that case,
coincides with all the field of algebraic functions
${\cal F}_n$.$\Box$

\section{Remarks concerning Theorem \ref{cor1.2}}
\label{sec+} 
\setcounter{equation}{0} The purpose of this
section is to show how the convergence result given by
Theorem
\ref{cor1.2} can be derived from the arguments of \S
\ref{sec5} more simply than the arguments given in \S
\ref{sec6}. Thus, let
$f:(M,0)\rightarrow (M',0)$ be a formal nondegenerate
CR map from a real analytic generic submanifold into a
real algebraic one. We also assume that $M'$ is
holomorphically nondegenerate. Then, by Proposition
\ref{propbis}, we know that for each component $f_j$
of $f$,
$j=1,\ldots,n$, there exists a positive integer $k_j$
and a family of convergent power series
$\delta_{i,j}=\delta_{i,j}\left((\Lambda_{\gamma})_{|\gamma|\leq
l(M')},z,w\right)\in
\C\{(\Lambda_{\gamma}-\partial^{\gamma}\bar{f}(0))_{|\gamma|\leq
l(M')},z,w\}$, $i=0,\ldots,k_j$, such that the formal
identity
\begin{equation}\label{eqRI}
\sum_{i=0}^{k_j}\delta_{i,j}\left((\partial^{\gamma}\bar{f}(w))_{|\gamma|\leq
l(M')},z,w\right) \left(f_j(z)\right)^i =0
\end{equation}
holds for $(z,w)\in {\cal M}$, with
$\delta_{k_j,j}\left((\partial^{\gamma}\bar{f}(w))_{|\gamma|\leq
l(M')},z,w\right)\not \equiv 0$ on ${\cal M}$. Here,
we recall that $l(M')$ is the Levi-type of $M'$ as
in
Remark \ref{rk1} and that ${\cal M}$ is the
complexification of
$M$ as defined in \S \ref{ssec2.1}.
Equation (\ref{eqRI}) means that for each
$j=1,\ldots,n$, $f_j(z)$ satisfies the statement (ii)
of Theorem \ref{th5.1}, with associated formal power
series mapping
$X(w)=(\partial^{\gamma}\bar{f}(w))_{|\gamma|\leq
l(M')}$. Thus, if we apply Proposition \ref{prop5.1}
to 
$h(z)=f_j(z)$ (and to
$X(w)=(\partial^{\gamma}\bar{f}(w))_{|\gamma|\leq
l(M')}$) for $j=1,\ldots,n$, we obtain the following
result.
\begin{proposition}\label{prop2+} Let $M\subset \C^n$ be a real analytic generic submanifold through the origin
and $M'\subset \C^n$ be a real algebraic generic 
submanifold through the origin (with the same CR
dimension).  Let
$f:M\rightarrow M'$ be a formal nondegenerate CR map
between $M$ and
$M'$ and assume that $M'$ is holomorphically
nondegenerate.  Then, for any multi-index $\mu\in
\N^n$ and for any $j \in
\{1,\ldots,n\}$, there exists a positive integer
$l(\mu,j)$, a family of convergent power series
$\delta_{i \mu j}=\delta_{i \mu
j}\left((\Lambda_{\gamma})_{|\gamma|\leq
l(M')+|\mu|},z,w\right)\in
\C\{(\Lambda_{\gamma}-\partial^{\gamma}
\bar{f}(0))_{|\gamma|\leq
l(M')+|\mu|},z,w\}$, $i=0,\ldots,l(\mu,j)$,  such
that the formal identity
$$\sum_{i=0}^{l(\mu,j)} \delta_{i \mu
j}\left((\partial^{\gamma}\bar{f}(w))_{|\gamma|\leq
l(M')+|\mu|},z,w\right)\left(\partial^{\mu}f_{j}(z)\right)^i=0,$$
holds for $(z,w)\in {\cal M}$ with $\delta_{l(\mu,j)
\mu j}\left((\partial^{\gamma}\bar{f}(w))
_{|\gamma|\leq l(M')+|\mu|},z,w\right)\not =0$ on
${\cal M}$.
\end{proposition} 
If furthermore $M$ is assumed to be minimal at 0,
then, in view of (\ref{eqRI}), we may apply
Proposition \ref{prop5.2} to $h(z)=f_j(z)$ and
$X(w)=(\partial^{\gamma}\bar{f}(w))_{|\gamma|\leq
l(M')}$, for $j=1,\ldots,n$. This gives the following
proposition.
\begin{proposition}\label{prop6.2} Let $M\subset \C^n$ be a real analytic generic
submanifold through the origin and $M'\subset \C^n$ be a real algebraic generic submanifold through the origin.
Let $f:M\rightarrow M'$ be a formal nondegenerate CR map between $M$ and $M'$ and 
assume that $M$ is minimal at 0 and
that $M'$ is holomorphically nondegenerate. Let $N$ be
the CR dimension of $M$ (and of $M'$) and $v_j$,
$j\in \N$, be the Segre sets mappings for $M$ as
defined by $(\ref{eqss1})$ and $(\ref{eqss2})$. Then,
for any multi-index
$\mu\in
\N^n$, for any
$d\in
\N$ and for any $j \in \{1,\ldots,n\}$, there exist
two positive integers $l=l(\mu,d,j),\
p=p(\mu,d,j)$, and a family of convergent power
series
$\psi_{\nu j}^{\mu
d}=\psi^{\mu
d}_{\nu j}\left((\Lambda_{\gamma})_{|\gamma|
\leq
p},z,w\right)\in
\C\{(\Lambda_{\gamma}-\partial^{\gamma}\bar{f}(0))
_{|\gamma|\leq
p},z,w\}$, $\nu=0,\ldots,l$, such that the formal
identity
 $$\sum_{\nu=0}^l
\psi_{\nu j}^{\mu d}\left(
((\partial^{\gamma}\bar{f})\circ
\bar{v}_d)_{|\gamma|\leq p},\bar{v}_d,v_{d+1}\right)
\left((\partial^{\mu}f_{j})\circ
v_{d+1}\right)^{\nu}\equiv 0$$ holds in the ring of
formal power series in $(d+1)N$ determinates and
such that
$\psi_{lj}^{\mu
d}\left(((\partial^{\gamma}\bar{f})\circ
\bar{v}_d)_{|\gamma|\leq
p},\bar{v}_d,v_{d+1}\right)\not
\equiv 0$.  \end{proposition}
From Proposition \ref{prop6.2}, one sees that 
the convergence of the mapping $f$ (under the
assumptions of Theorem
\ref{cor1.2}) follows from successive applications of
Lemma
\ref{lem5.4}. Indeed, for $d=0$, Proposition
\ref{prop6.2}, and Lemma \ref{lem5.4} yield the
convergence of
$f$ and of all its jets on the first Segre set. From
this, Proposition
\ref{prop6.2} and Lemma \ref{lem5.4}, we obtain the
convergence of $f$ and of all its jets on the
second Segre set, and so on. This leads to the
convergence of $h$ on the $d_0$-th Segre set, where
$d_0$ is given by Theorem \ref{th2.1}. As in the 
proof of
Theorem \ref{th5.1}, this
implies the convergence of $f$ under the assumptions
 of Theorem \ref{cor1.2}.

\section{Partial convergence of formal CR maps}
\label{secpc}
In this section, as in \cite{M3}, we indicate
several results which show how Theorem \ref{th1.1}
can be seen as a result of partial convergence for
formal nondegenerate CR maps. Before explaining
what we mean by this, we have to recall several facts.

Let
$M$ be a real analytic generic submanifold in $\C^n$
and $p\in M$. Let $\K(p)$ be the quotient field of
$\C\{z-p\}$, and
$H(M,p)$ be the vector space over $\K(p)$
consisting of the germs at $p$ of (1,0) vector
fields, with meromorphic coefficients, tangent to
$M$ (near $p$). The {\it degeneracy}
of
$M$ at
$p$, denoted $d(M,p)$, is defined to be the dimension
of
$H(M,p)$ over $\K(p)$. It is shown in \cite{BR, BER1}
that the mapping $M\ni p \mapsto d(M,p)\in
\{0,\ldots,n\}$ is constant on any connected
component of $M$. Consequently, if $M$ is a connected
real analytic generic submanifold, one can define its
degeneracy $d(M)$ to be the degeneracy $d(M,q)$ at
any point $q\in M$.  Observe that the germ 
$(M,p)$, $p\in M$, is holomorphically nondegenerate
if and only if
$d(M)=d(M,p)=0$.

If $f$ is a formal nondegenerate CR map as in
Theorem \ref{th1.1}, $f$ can or cannot be
convergent. The following result, which is of
more interest when $f$ is not convergent, shows
however that the map $f$ is partially convergent in
the following sense.

\begin{theorem}\label{thpc0}
Let $f:(M,0)\rightarrow (M',0)$ be a formal
nondegenerate CR map between two germs at 0 in $\C^n$
of real analytic generic submanifolds. Assume that $M$
is minimal at 0 and that $M'$ is real algebraic. Let
$d(M')$ be the degeneracy of $M'$. Then, there exists
a holomorphic (algebraic) mapping (independent
of $f$) ${\cal G}(\omega)=({\cal
G}_1(\omega),\ldots,{\cal G}_{n-d(M')}(\omega))$
defined near $0\in \C^n$ of generic rank $n-d(M')$
such that ${\cal G}\circ f$ is
convergent. 
\end{theorem}
{\it Proof.}
We use again the notations of \S \ref{sec1} and \S
\ref{ssec4.1}. As in the proof of Theorem
\ref{cor1.2}, by using the expansion (\ref{taylor}),
we obtain for
$\nu =1,\ldots,c$,
$$\bar \Phi_{\nu}'(f(z),\theta)=\sum_{\beta \in
\N^N}(q_{\beta,\nu}\circ f)(z)\, \theta^{\beta}.$$
Recall also that the $q_{\beta,\nu}(\omega)$ are
algebraic functions. By Theorem \ref{th1.1}, we have
that for any $\beta \in \N^N$ and for any
$\nu=1,\ldots,c$, $q_{\beta, \nu}\circ f$ is
convergent in some neighborhood of $0$ in
$\C^n$. According to \cite{BR, BER1}, we can choose
$q_{\beta^1,\nu^1}(\omega),\ldots,q_{\beta^r,\nu^r}(\omega)$,
$r=n-d(M')$, of generic rank $n-d(M')$ in a
neighborhood of 0 in $\C^n$. Then, if we define
${\cal G}_i(\omega)=q_{\beta^i,\nu^i}(\omega)$, for
$i=1,\ldots,n-d(M')$, we obtain the desired statement
of the Theorem.$\Box$

A suitable invariant which also measures the lack of
convergence of a given formal (holomorphic) mapping is
its so-called {\it transcendence degree}. We recall
first how such an invariant is defined ({\it
cf.} \cite{M3}). Let ${\cal
H}:(\C^N,0)\rightarrow (\C^{N'},0)$ be formal holomorphic
mapping, $N,N'\geq 1$, and 
$V$ be a complex analytic set through the origin in $\C^N\times
\C^{N'}$. Assume that $V$ is given near the origin in $\C^{N+N'}$
by $V=\{(x,y)\in \C^N\times
\C^{N'}:{\varrho}_1(x,y)=\ldots={\varrho}_q(x,y)=0\}$,
${\varrho}_i(x,y)\in \C\{x,y\}$, $i=1,\ldots,q$. Then,
the graph of
${\cal H}$ is said to be {\it formally contained} in $V$
if ${\varrho}_1(x,{\cal
H}(x))=\ldots={\varrho}_q(x,{\cal H}(x))=0$ in
$\C[[x]]$. Furthermore, if
$V_{\cal H}$ is the germ of the complex analytic set
through the origin in $\C^{N+N'}$ defined as the
intersection of all the complex analytic sets
through the origin in $\C^{N+N'}$ formally containing
the graph of
${\cal H}$, then the {\it transcendence degree} of
${\cal H}$ is defined to be the integer ${\rm
dim}_{\C}\ V_{\cal H}-N$. This definition is in
complete analogy with the one introduced in
\cite{CPS2} in the ${\cal C}^{\infty}$ mapping
problem. The following proposition from \cite{M3}
shows the relevance of the previous concept and why
this transcendence degree is related to the
convergence properties of formal mappings.

\begin{proposition}\label{propc} Let ${\cal
H}:(\C^N,0)\rightarrow (\C^{N'},0)$ be formal
holomorphic mapping. Then, the following conditions
are equivalent:\\ i)
${\cal H}$ is convergent.\\ ii) The transcendence
degree of ${\cal H}$ is zero.  \end{proposition}  The
following is a consequence of Theorem
\ref{thpc0}.

\begin{corollaire}\label{thpc} Let
$f:(M,0)\rightarrow (M',0)$ be a formal CR mapping
between two germs at 0 of real analytic generic
submanifolds.  Assume that
$M$ is minimal at 0, $M'$ is real algebraic and that
$f$ is nondegenerate, i.e.  $J_f\not \equiv 0$. 
Denote by ${\cal D}_f$ the transcendence degree of
$f$. Then, ${\cal D}_f\leq d(M')$, where $d(M')$
denotes the degeneracy of $M'$. Equivalently, there
exists a complex analytic set of (pure) dimension
$n+d(M')$ which contains formally the graph of $f$. 
\end{corollaire}

{\it Proof.} The proof is similar to the one
given in
\cite{M3}. For the sake of completeness, we recall it.

From Theorem \ref{thpc0}, there exists ${\cal
G}(\omega)=({\cal G}_1(\omega),\ldots,{\cal
G}_{n-d(M')}(\omega))\in (\C\{\omega\})^{n-d(M')}$ of
generic rank $n-d(M')$ such that ${\cal G}\circ f$ is
convergent. Put $\delta_j(z):=({\cal G}_j\circ
f)(z)\in \C\{z\}$, for $j=1,\ldots,n-d(M')$. Then, the
graph of
$f$ is formally contained in the complex analytic set
$$A=\{(z,\omega)\in
(\C^{2n},0):{\cal
G}_{j}(\omega)=\delta_{j}(z),\
j=1,\ldots,n-d(M')\}.$$ Let
$A=\cup_{i=1}^k\Gamma_i$ be the decomposition of
$A$ into irreducible components.  For any positive integer
$\sigma$, one can find, according to the Artin
approximation theorem \cite{A}, a convergent power
series mapping
$f^{\sigma}(z)\in (\C\{z\})^n$ defined in some small
neighborhood
$U^{\sigma}$ of $0$ in $\C^n$, which agrees with $f(z)$
up to order
$\sigma$ (at 0) and such that the graph of
$f^{\sigma}$, denoted
$G(f^{\sigma})$, is contained in $A$.  Since
$G(f^{\sigma})$ is contained in $A$, it must be
contained in an irreducible component of $A$.  Thus,
at least one subsequence of
$(f^{\sigma})_{\sigma \in \N^*}$ is contained in one of
such irreducible components, say $\Gamma_1$.  We can
assume without loss of generality that such a
subsequence is
$(f^{\sigma})_{\sigma \in \N^*}$ itself.  We first
observe that this implies that the graph of $f$ is
formally contained in
$\Gamma_1$.  Moreover, since
$f$ is a formal nondegenerate map, for $\sigma_0$ large
enough,  the family $(f^{\sigma})_{\sigma \geq
\sigma_0}$ is also a family of holomorphic maps of
generic rank
$n$.  In particular, this implies that the generic
rank of the family of holomorphic functions
$$\left(({\cal G}_{i}\circ
f^{\sigma_0})(z)\right)_{1\leq i\leq n-d(M')},$$ is
$n-d(M')$.  As a consequence, if $z_0$ is close enough
to 0 in
$\C^n$ and is chosen so that the rank of the preceding
family at $z_0$ equals
$n-d(M')$, the implicit function theorem shows that $A$
is an
$n+d(M')$-dimensional complex submanifold near
$(z_0,f^{\sigma_0}(z_0))\in
\Gamma_1$.  Since $\Gamma_1$ is irreducible, it is
pure-dimensional, and thus,
$\Gamma_1$ is an $n+d(M')$ pure-dimensional complex analytic set
containing formally the graph of $f$.  By definition of the
transcendence degree, this implies that ${\cal D}_f\leq
d(M')$.$\Box$

\begin{remark}\label{rk11}
One should observe that
Theorem
\ref{cor1.2} also follows from Corollary \ref{thpc}.
Indeed, if, in Corollary
\ref{thpc}, $M'$ is furthermore assumed to be
holomorphically nondegenerate, then $d(M')=0$ and
thus the transcendence degree of $f$ is zero. By
Proposition \ref{propc}, this implies that $f$ is
convergent.
\end{remark}
 
\section{Concluding remarks}\label{sec7} In this last section, we indicate how Theorem \ref{th5.1} can be applied to the study
of the convergence of formal mappings between real
analytic CR manifolds in complex spaces of possibly
different dimensions.  Our last result will be
expressed by means of a standard nondegeneracy
condition which takes its source in
\cite{P, W1, DF}. The situation is the following one.

Let $f :  (M,0)\rightarrow (M',0)$ be a
formal CR mapping between two germs at 0 of real
analytic generic submanifolds in $\C^n$ and $\C^{n'}$
respectively, $n,n'\geq 2$. (We wish to mention that
all the following considerations are also valid for a
target real analytic set, but for simplicity, we
restrict our attention to generic manifolds.)  We
shall use the notations defined in \S
\ref{sec2} for $M$. In particular, the CR dimension of
$M$ is $N$ and its real codimension is $c$. Following
the terminology of
\cite{CPS2}, we define the {\it characteristic variety
of $f$} at $0\in \C^{n'}$ as follows. If $M'$ is a real
analytic generic submanifold through
$0$ as above, of CR dimension $N'$ and of real
codimension
$c'$, we can assume that it is given near $0$ by
$$M'=\{\zeta
\in (\C^{n'},0):\rho' (\zeta,\bar{\zeta})=0\},$$ 
with $$\bar{\partial} \rho'_1 \wedge \ldots \wedge
\bar{\partial} \rho'_{c'} \not =0,\ {\rm on}\ M'.$$
Here, $\rho'=(\rho'_1,\ldots,\rho'_{c'})$ is a set
of local real analytic defining functions for $M'$ near
$0\in \C^{n'}$. Consider the vector fields ${\cal
L}_j$, $j=1,\ldots,N$, tangent to the complexification
${\cal M}$ of the source manifold $M$ as defined by
(\ref{vf2}).  It will be better to see these vector
fields, for fixed $z\in (\C^n,0)$, as a basis of
holomorphic vector fields tangent to the Segre variety
$Q_{\bar{z}}$, and thus, we shall denote them ${\cal
L}_{(z,w)}^j$ for
$j=1,\ldots,N$. For any multi-index
$\gamma =(\gamma_1,\ldots,\gamma_N) \in
\N^N$, we define ${\cal L}_{(z,w)}^{\gamma}=\left[{\cal
L}_{(z,w)}^{1}\right]^{\gamma_1}\ldots \left[{\cal
L}_{(z,w)}^{N}\right]^{\gamma_N}$.  Finally, for
any multi-index $\gamma \in \N^N$, let $\Xi_{\gamma}$
be the
$\C^{c'}$ formal map defined by
\begin{equation}\label{eqchi}
\C^n \times \C^{n}\times
\C^{n'}\ni (z,w,\zeta)\mapsto {\cal L}_{(z,w)}^
{\gamma}\rho'(\zeta,\bar{f}(w))\in \C^{c'}.
\end{equation}  Observe that there
exists a $\C^{c'}$-valued  convergent
power series mapping ${\cal N}_{\gamma}={\cal N}
_{\gamma}\left((\Lambda_{\beta})_{|\beta|\leq
|\gamma|},z,w,\zeta \right)\in
\left(\C\{\Lambda_0,z,w,\zeta\}[(\Lambda_{\beta}))_{1\leq
|\beta|\leq
|\gamma|}]\right)^{c'}$ such
that
\begin{equation}\label{eqchibis}
\Xi_{\gamma}(z,w,\zeta)={\cal
N}_{\gamma}\left((\partial^{\beta}
\bar{f}(w))_{|\beta|\leq |\gamma|},z,w,\zeta\right),\
{\rm in}\ \C[[z,w,\zeta]].
\end{equation}
The characteristic variety
of
$f$ at $0\in \C^n$ is then defined to be the germ at 
$0\in \C^{n'}$ of
the complex analytic set $${\cal C}(f,0)=\{\zeta\in
(\C^{n'},0):\Xi_{\gamma}(0,0,\zeta)=0,\gamma
\in \N^N\}.$$ 
This set is the infinitesimal analog of the usual
determinacy set for holomorphic mappings 
$$\{\zeta \in (\C^{n'},0):f(Q_0)\subset
Q'_{\zeta}\},$$ 
where $Q'_{\zeta}$ is the Segre variety associated
to $M'$ and $\zeta \in (\C^{n'},0)$. With this in mind,
we have the following natural result. 
\begin{theorem}\label{th7.1} Let $f :  (M,0)
\rightarrow (M',0)$ be a formal
mapping between two germs at 0 of real analytic
generic  submanifolds in $\C^n$ and $\C^{n'}$
respectively. Assume that
$M$ is minimal at
$0$ and that the characteristic variety ${\cal
C}(f,0)$ is zero-dimensional.  Then $f$ is
convergent.  \end{theorem} 
{\it Proof.} Since $f$ maps formally
$M$ to $M'$, we have
\begin{equation}\label{end1}
\rho'(f(z),\bar{f}(w))=0,
\end{equation}
as a formal power series identity for $(z,w)\in
{\cal M}$. Thus, if, for $\gamma \in \N^N$, we apply
${\cal L}_{(z,w)}^{\gamma}$ to (\ref{end1}), it follows
from the definition of $\Xi_{\gamma}$ given in
(\ref{eqchi}), that
\begin{equation}\label{end2}
{\cal
L}^{\gamma}_{(z,w)}\left(\rho'(f(z),\bar{f}(w))
\right)=\Xi_{\gamma}(z,w,f(z))=0,
\end{equation} 
for $(z,w)\in {\cal M}$. Observe that it follows from
(\ref{eqchibis}) and (\ref{end2}) that
\begin{equation}\label{eq7.1}
{\cal N}_{\gamma}
\left((\partial^{\beta} \bar{f}(w))_{|\beta|\leq
|\gamma|},z,w,f(z)\right)=0,\ {\rm for}\ (z,w)\in {\cal
M}.
 \end{equation} 
Since the
characteristic variety ${\cal C}(f,0)$ is
zero-dimensional, in view of (\ref{eqchibis}), 
the holomorphic mapping ${\cal I}_k$
(\footnote{$d(i,j)=
 {\rm Card}\{\alpha \in
\N^{i}:|\beta|\leq j\}$, $i,j\in \N^*$.})
\begin{multline}\label{end4}
\C^{n'd(n,k)+2n+n'} \ni
((\Lambda_{\beta})_{|\beta|\leq k},z,w,\zeta)
\mapsto \\
\left((\Lambda_{\beta})_{|\beta|\leq k},z,w,\left({\cal
N}_{\gamma}((\Lambda_{\beta})_{|\beta|\leq
|\gamma|},z,w,\zeta)\right)_{|\gamma|\leq k}\right)
\end{multline} is
finite-to-one near
$J_0=((\partial^{\beta}\bar{f}(0))_{|\beta|\leq
k},0,0,0)\in \C^{n'd(n,k)+2n+n'}$ for $k$ large enough.
It then follows from \cite{Gu} (p.15) (see also
\cite{MT}) that, for any $j=1,\ldots,n'$, $\zeta_j$ is
integral over the {\it ring} formed by all the
convergent power series of the form $$({\cal B}\circ
{\cal I}_k)\left((\Lambda_{\beta})_{|\beta|\leq
k},z,w,\zeta\right),$$
${\cal B}$ running over all the convergent power
series centered at
$$J_1=((\partial^{\beta}\bar{f}(0))_{|\beta|\leq
k},0,0,0)\in \C^{n'd(n,k)+2n+c'd(N,k)}.$$
Explicitly,
for any
$j=1,\ldots,n'$, there exists a positive integer
$\nu_j$ and convergent power series ${\cal B}_{t,j}$
near $J_1$, $t=0,\ldots,\nu_j-1$,
such that the following identities hold in a
neigborhood of $J_0$: 
\begin{equation}\label{eq7.2}
\zeta_j^{\nu_j}+\sum_{t<\nu_j}({\cal B}_{t,j}\circ
{\cal I}_k)\left((\Lambda_{\beta})_{|\beta|\leq
k},z,w,\zeta \right)\ \zeta_j^{t}=0.  \end{equation} 
Putting
$\zeta =f(z)$ and $(\Lambda_{\beta})_{|\beta|\leq
k}=(\partial^{\beta}\bar{f}(w))_{|\beta|\leq k}$ for
$(z,w)\in (\C^{2n},0)$ in (\ref{eq7.2}), we obtain that
for $(z,w)\in (\C^{2n},0)$ and for $j=1,\ldots,n'$, the
following formal identities hold in $\C[[z,w]]$:
\begin{equation}\label{end6}
\left(f_j(z)\right)^{\nu_j}+\sum_{t<\nu_j}({\cal
B}_{t,j}\circ {\cal
I}_k)\left((\partial^{\beta}\bar{f}(w))_{|\beta|\leq
k},z,w,f(z)\right)\ \left(f_j(z)\right)^{t}=0.
\end{equation}
But in view of the definition of
${\cal I}_k$ given in (\ref{end4}) and in view of
(\ref{eq7.1}),  we have for $j=1,\ldots,n'$, for 
$t=0,\ldots,\nu_j-1$ and for $(z,w)\in {\cal M}$,
\begin{multline}
({\cal
B}_{t,j}\circ {\cal
I}_k)\left((\partial^{\beta}\bar{f}(w))_{|\beta|\leq
k},z,w,f(z)\right)=\\
 {\cal B}_{t,j}\left(
(\partial^{\beta}\bar{f}(w))_{|\beta|\leq
k},z,w,\left({\cal
N}_{\gamma}((\partial^{\beta}\bar{f}(w))_{|\beta|\leq
|\gamma|},z,w,f(z))\right)_{|\gamma|\leq k}\right)=\\
{\cal B}_{t,j}\left(
(\partial^{\beta}\bar{f}(w))_{|\beta|\leq
k},z,w,0\right).
\end{multline}
Thus, from (\ref{end6}), we obtain for $j=1,\ldots,n'$,
\begin{equation}\label{end8}
\left(f_j(z)\right)^{\nu_j}+
\sum_{t<\nu_j} {\cal B}_{t,j}\left(
(\partial^{\beta}\bar{f}(w))_{|\beta|\leq
k},z,w,0\right) \left(f_j(z)\right)^{t}=0,
\end{equation}
on ${\cal M}$. As a consequence, we see that for each
$j=1,\ldots,n'$, the formal holomorphic power series
$f_j(z)$, $j=1,\ldots,n'$, satisfies (ii) of Theorem
\ref{th5.1}. Since $M$ is minimal at 0, from that
theorem, we conclude that $f$ is convergent.$\Box$

\begin{remark}  It should be mentioned
that Theorem
\ref{th7.1} above could also be derived from the
techniques of \cite{BER4}.
\end{remark}

We conclude by mentioning several situations where
Theorem
\ref{th7.1} applies.  It contains the cases of formal invertible mappings of Levi-nondegenerate real analytic
hypersurfaces, finite mappings of minimal essentially 
finite real analytic generic manifolds or, more
generally, mappings with injective Segre homomorphim
(in the sense of \cite{BER4}) from minimal real
analytic generic manifolds into real analytic
essentially finite ones (the proof is contained in
\cite{BER4}).

\end{document}